%

\magnification\magstep1
\baselineskip15pt
\parskip 0.0pt plus 2.5pt

\newread\AUX\immediate\openin\AUX=\jobname.aux
\newcount\relFnno
\def\ref#1{\expandafter\edef\csname#1\endcsname}
\ifeof\AUX\immediate\write16{\jobname.aux gibt es nicht!}\else
\input \jobname.aux
\fi\immediate\closein\AUX



\def\ignore{\bgroup
\catcode`\;=0\catcode`\^^I=14\catcode`\^^J=14\catcode`\^^M=14
\catcode`\ =14\catcode`\!=14\catcode`\"=14\catcode`\#=14\catcode`\$=14
\catcode`\&=14\catcode`\'=14\catcode`\(=14\catcode`\)=14\catcode`\*=14
\catcode`+=14\catcode`\,=14\catcode`\-=14\catcode`\.=14\catcode`\/=14
\catcode`\0=14\catcode`\1=14\catcode`\2=14\catcode`\3=14\catcode`\4=14
\catcode`\5=14\catcode`\6=14\catcode`\7=14\catcode`\8=14\catcode`\9=14
\catcode`\:=14\catcode`\<=14\catcode`\==14\catcode`\>=14\catcode`\?=14
\catcode`\@=14\catcode`\A=14\catcode`\B=14\catcode`\C=14\catcode`\D=14
\catcode`\E=14\catcode`\F=14\catcode`\G=14\catcode`\H=14\catcode`\I=14
\catcode`\J=14\catcode`\K=14\catcode`\L=14\catcode`\M=14\catcode`\N=14
\catcode`\O=14\catcode`\P=14\catcode`\Q=14\catcode`\R=14\catcode`\S=14
\catcode`\T=14\catcode`\U=14\catcode`\V=14\catcode`\W=14\catcode`\X=14
\catcode`\Y=14\catcode`\Z=14\catcode`\[=14\catcode`\\=14\catcode`\]=14
\catcode`\^=14\catcode`\_=14\catcode`\`=14\catcode`\a=14\catcode`\b=14
\catcode`\c=14\catcode`\d=14\catcode`\e=14\catcode`\f=14\catcode`\g=14
\catcode`\h=14\catcode`\i=14\catcode`\j=14\catcode`\k=14\catcode`\l=14
\catcode`\m=14\catcode`\n=14\catcode`\o=14\catcode`\p=14\catcode`\q=14
\catcode`\r=14\catcode`\s=14\catcode`\t=14\catcode`\u=14\catcode`\v=14
\catcode`\w=14\catcode`\x=14\catcode`\y=14\catcode`\z=14\catcode`\{=14
\catcode`\|=14\catcode`\}=14\catcode`\~=14\catcode`\^^?=14
\Ignoriere}
\def\Ignoriere#1\;{\egroup}

\newcount\itemcount
\def\resetitem{\global\itemcount0}\resetitem
\newcount\Itemcount
\Itemcount0
\newcount\maxItemcount
\maxItemcount=0

\def\FILTER\fam\itfam\tenit#1){#1}

\def\Item#1{\global\advance\itemcount1
\edef\TEXT{{\it\romannumeral\itemcount)}}%
\ifx?#1?\relax\else
\ifnum#1>\maxItemcount\global\maxItemcount=#1\fi
\expandafter\ifx\csname I#1\endcsname\relax\else
\edef\testA{\csname I#1\endcsname}
\expandafter\expandafter\def\expandafter\testB\testA
\edef\testC{\expandafter\FILTER\testB}
\edef\testD{\csname0\testC0\endcsname}\fi
\edef\testE{\csname0\romannumeral\itemcount0\endcsname}
\ifx\testD\testE\relax\else
\immediate\write16{I#1 hat sich geaendert!}\fi
\expandwrite\AUX{\neverexpand\ref{I#1}{\TEXT}}\fi
\item{\ifx?#1?\relax\else\marginnote{I#1}\fi\TEXT}}

\def\today{\number\day.~\ifcase\month\or
  Januar\or Februar\or M{\"a}rz\or April\or Mai\or Juni\or
  Juli\or August\or September\or Oktober\or November\or Dezember\fi
  \space\number\year}
\font\sevenex=cmex7
\scriptfont3=\sevenex
\font\fiveex=cmex10 scaled 500
\scriptscriptfont3=\fiveex
\def\A{{\bf A}}
\def\G{{\bf G}}
\def\P{{\bf P}}

\def\phi{\varphi}
\def\epsilon{\varepsilon}
\def\theta{\vartheta}
\def\uauf{\lower1.7pt\hbox to 3pt{%
\vbox{\offinterlineskip
\hbox{\vbox to 8.5pt{\leaders\vrule width0.2pt\vfill}%
\kern-.3pt\hbox{\lams\char"76}\kern-0.3pt%
$\raise1pt\hbox{\lams\char"76}$}}\hfil}}

\font\BF=cmbx10 scaled \magstep2

\def\title#1{\par
{\baselineskip1.5\baselineskip\rightskip0pt plus 5truecm
\leavevmode\vskip0truecm\noindent\BF #1\par}
\vskip1truecm
\leftline{\font\CSC=cmcsc10
{\CSC Friedrich Knop}}
\leftline{Department of Mathematics, Rutgers University, Piscataway NJ
08854, USA}
\leftline{knop@math.rutgers.edu}
\vskip1truecm
\par}

\def\cite#1{\expandafter\ifx\csname#1\endcsname\relax
{\bf?}\immediate\write16{#1 ist nicht definiert!}\else\csname#1\endcsname\fi}
\def\expandwrite#1#2{\edef\next{\write#1{#2}}\next}
\def\neverexpand{\noexpand\noexpand\noexpand}
\def\strip#1\ {}
\def\ncite#1{\expandafter\ifx\csname#1\endcsname\relax
{\bf?}\immediate\write16{#1 ist nicht definiert!}\else
\expandafter\expandafter\expandafter\strip\csname#1\endcsname\fi}
\newwrite\AUX
\immediate\openout\AUX=\jobname.aux
\font\eightrm=cmr8\font\sixrm=cmr6
\font\eighti=cmmi8
\font\eightit=cmti8
\font\eightbf=cmbx8
\font\eightcsc=cmcsc10 scaled 833
\def\eightpoint{%
\textfont0=\eightrm\scriptfont0=\sixrm\def\rm{\fam0\eightrm}%
\textfont1=\eighti
\textfont\bffam=\eightbf\def\bf{\fam\bffam\eightbf}%
\textfont\itfam=\eightit\def\it{\fam\itfam\eightit}%
\def\csc{\eightcsc}%
\setbox\strutbox=\hbox{\vrule height7pt depth2pt width0pt}%
\normalbaselineskip=0,8\normalbaselineskip\normalbaselines\rm}
\newcount\absFnno\absFnno1
\write\AUX{\relFnno1}
\newif\ifMARKE\MARKEtrue
{\catcode`\@=11
\gdef\footnote{\ifMARKE\edef\@sf{\spacefactor\the\spacefactor}\/%
$^{\cite{Fn\the\absFnno}}$\@sf\fi
\MARKEtrue
\insert\footins\bgroup\eightpoint
\interlinepenalty100\let\par=\endgraf
\leftskip=0pt\rightskip=0pt
\splittopskip=10pt plus 1pt minus 1pt \floatingpenalty=20000\smallskip
\item{$^{\cite{Fn\the\absFnno}}$}%
\expandwrite\AUX{\neverexpand\ref{Fn\the\absFnno}{\neverexpand\the\relFnno}}%
\global\advance\absFnno1\write\AUX{\advance\relFnno1}%
\bgroup\strut\aftergroup\@foot\let\next}}
\skip\footins=12pt plus 2pt minus 4pt
\dimen\footins=30pc
\output={\plainoutput\immediate\write\AUX{\relFnno1}}
\newcount\Abschnitt\Abschnitt0
\def\beginsection#1. #2 \par{%
\advance\Abschnitt1
\edef\Aname{\number\Abschnitt}
\vskip0pt plus.10\vsize\penalty-250
\vskip0pt plus-.10\vsize\bigskip\vskip\parskip
\expandafter\ifx\csname#1\endcsname\Aname\relax\else
\immediate\write16{#1 hat sich geaendert!}\fi
\expandwrite\AUX{\neverexpand\ref{#1}{\Aname}}
\leftline{\marginnote{#1}\bf\Aname. \ignorespaces#2}%
\nobreak\smallskip\noindent\SATZ1\GNo0%
\write\TOC{\string\toc{\Aname. #2}{\the\pageno}}}
\def\References{%
\vskip0pt plus.10\vsize\penalty-250
\vskip0pt plus-.10\vsize\bigskip\vskip\parskip
\leftline{\bf References}\nobreak\smallskip\noindent%
\write\TOC{\string\toc{References}{\the\pageno}}}
\newcount\appcount

\def\Appendix#1. #2\par{%
\appcount0
\vskip0pt plus.10\vsize\penalty-250
\vskip0pt plus-.10\vsize\bigskip\bigskip\bigskip\vskip\parskip
\xdef\AppName{#1}
\noindent{\BF Appendix #1: #2}\par
\nobreak\medskip\noindent
\write\TOC{\string\toc{Appendix #1: #2}{\the\pageno}}}

\def\appsection#1. #2 \par{%
\advance\appcount1
\edef\Aname{\AppName\number\appcount}
\vskip0pt plus.10\vsize\penalty-250
\vskip0pt plus-.10\vsize\bigskip\vskip\parskip
\expandwrite\AUX{\neverexpand\ref{#1}{A}}
\leftline{\marginnote{#1}\bf\Aname. \ignorespaces#2}%
\nobreak\smallskip\noindent\SATZ1\GNo0%
\write\TOC{\string\toc{\hskip20pt\Aname. #2}{\the\pageno}}}
\newif\ifmarginalnotes\marginalnotesfalse
\newif\ifmarginalwarnings\marginalwarningstrue

\def\marginnote#1{\ifmarginalnotes\hbox to 0pt{\eightpoint\hss #1\ }\fi}

\def\strutdepth{\dp\strutbox}
\def\Randbem#1#2{\ifmarginalwarnings
{#1}\strut
\setbox0=\vtop{\eightpoint
\rightskip=0pt plus 6mm\hfuzz=3pt\hsize=16mm\noindent\leavevmode#2}%
\vadjust{\kern-\strutdepth
\vtop to \strutdepth{\kern-\ht0
\hbox to \hsize{\kern-16mm\kern-6pt\box0\kern6pt\hfill}\vss}}\fi}

\def\Zitat!{\Randbem{\bf?}{\bf Zitat}}

\newcount\SATZ\SATZ1
\def\proclaim #1. #2\par{\ifdim\lastskip<\medskipamount\removelastskip
\medskip\fi
\noindent{\bf#1.\ }{\it#2\Par}
\ifdim\lastskip<\medskipamount\removelastskip\goodbreak\medskip\fi}
\def\rproclaim #1. #2\par{\ifdim\lastskip<\medskipamount\removelastskip
\medskip\fi
\noindent{\bf#1.\ }{\rm#2\Par}
\ifdim\lastskip<\medskipamount\removelastskip\goodbreak\medskip\fi}
\def\Aussage#1{\expandafter\def\csname#1\endcsname##1.{\resetitem
\ifx?##1?\relax\else
\edef\TEST{#1\penalty10000\ \Aname.\number\SATZ}
\expandafter\ifx\csname##1\endcsname\TEST\relax\else
\immediate\write16{##1 hat sich geaendert!}\fi
\expandwrite\AUX{\neverexpand\ref{##1}{\TEST}}\fi
\proclaim {\marginnote{##1}\Aname.\number\SATZ. #1\global\advance\SATZ1}.}}
\def\rAussage#1{\expandafter\def\csname#1\endcsname##1.{\resetitem
\ifx?##1?\relax\else
\edef\TEST{#1\penalty10000\ \Aname.\number\SATZ}
\expandafter\ifx\csname##1\endcsname\TEST\relax\else
\immediate\write16{##1 hat sich geaendert!}\fi
\expandwrite\AUX{\neverexpand\ref{##1}{\TEST}}\fi
\rproclaim {\marginnote{##1}\Aname.\number\SATZ. #1\global\advance\SATZ1}.}}
\Aussage{Theorem}
\Aussage{Proposition}
\Aussage{Corollary}
\Aussage{Lemma}
\def\Proof:{\par\noindent{\it Proof:}}
\rAussage{Definition}
\def\Remark:{\ifdim\lastskip<\medskipamount\removelastskip\medskip\fi
\noindent{\bf Remark:}}
\def\Remarks:{\ifdim\lastskip<\medskipamount\removelastskip\medskip\fi
\noindent{\bf Remarks:}}
\def\Example:{\ifdim\lastskip<\medskipamount\removelastskip\medskip\fi
\noindent{\bf Example:}}
\def\Examples:{\ifdim\lastskip<\medskipamount\removelastskip\medskip\fi
\noindent{\bf Examples:}}
\font\la=lasy10
\def\strich{\hbox{$\vcenter{\hbox
to 1pt{\leaders\hrule height -0,2pt depth 0,6pt\hfil}}$}}
\def\dashedrightarrow{\hbox{%
\hbox to 0,5cm{\leaders\hbox to 2pt{\hfil\strich\hfil}\hfil}%
\kern-2pt\hbox{\la\char\string"29}}}

\def\Bindestrich{\penalty10000-\hskip0pt}
\let\_=\Bindestrich
\def\.{{\sfcode`.=1000.}}

\def\Par{\par}
\def\:={\mathrel{\raise0,9pt\hbox{.}\kern-2,77779pt
\raise3pt\hbox{.}\kern-2,5pt=}}
\def\=:{\mathrel{=\kern-2,5pt\raise0,9pt\hbox{.}\kern-2,77779pt
\raise3pt\hbox{.}}} \def\mod{/\mskip-5mu/}
\def\into{\hookrightarrow}
\def\pfeil{\rightarrow}

\def\Pfeil{\longrightarrow}
\def\pf#1{\buildrel#1\over\rightarrow}
\def\Pf#1{\buildrel#1\over\longrightarrow}

\def\Ugleich{\hbox{$\cup$\kern.5pt\vrule depth -0.5pt}}
\def\|#1|{\mathop{\rm#1}\nolimits}
\def\<{\langle}
\def\>{\rangle}
\let\Times=\times
\def\times{\mathop{\Times}}
\let\Otimes=\otimes
\def\otimes{\mathop{\Otimes}}
\catcode`\@=11
\def\hex#1{\ifcase#1 0\or1\or2\or3\or4\or5\or6\or7\or8\or9\or A\or B\or
C\or D\or E\or F\else\message{Warnung: Setze hex#1=0}0\fi}
\def\fontdef#1:#2,#3,#4.{%
\alloc@8\fam\chardef\sixt@@n\FAM
\ifx!#2!\else\expandafter\font\csname text#1\endcsname=#2
\textfont\the\FAM=\csname text#1\endcsname\fi
\ifx!#3!\else\expandafter\font\csname script#1\endcsname=#3
\scriptfont\the\FAM=\csname script#1\endcsname\fi
\ifx!#4!\else\expandafter\font\csname scriptscript#1\endcsname=#4
\scriptscriptfont\the\FAM=\csname scriptscript#1\endcsname\fi
\expandafter\edef\csname #1\endcsname{\fam\the\FAM\csname text#1\endcsname}
\expandafter\edef\csname hex#1fam\endcsname{\hex\FAM}}
\catcode`\@=12 

\fontdef Ss:cmss10,cmss7,.
\fontdef Fr:eufm10,eufm7,eufm5.

\def\fE{{\Fr E}}


\def\fP{{\Fr P}}
\def\fQ{{\Fr Q}}

\def\fS{{\Fr S}}
\def\fT{{\Fr T}}

\fontdef bbb:msbm10,msbm7,msbm5.
\fontdef mbf:cmmib10,cmmib7,.
\fontdef msa:msam10,msam7,msam5.
\def\CC{{\bbb C}}
\def\FF{{\bbb F}}

\def\NN{{\bbb N}}
\def\QQ{{\bbb Q}}

\def\ZZ{{\bbb Z}}
\def\cA{{\cal A}}\def\cB{{\cal B}}\def\cC{{\cal C}}\def\cD{{\cal D}}

\def\cI{{\cal I}}
\def\cN{{\cal N}}\def\cO{{\cal O}}
\def\cT{{\cal T}}

\mathchardef\leer=\string"0\hexbbbfam3F
\mathchardef\subsetneq=\string"3\hexbbbfam24
\mathchardef\semidir=\string"2\hexbbbfam6E
\mathchardef\dirsemi=\string"2\hexbbbfam6F
\mathchardef\haken=\string"2\hexmsafam78
\mathchardef\auf=\string"3\hexmsafam10
\let\OL=\overline
\def\overline#1{{\hskip1pt\OL{\hskip-1pt#1\hskip-.3pt}\hskip.3pt}}

\def\eq{{\overline{e}}}
\def\fQ{{\overline{f}}}

\def\gq{{\overline{g}}}

\def\mq{{\overline{m}}}

\def\rq{{\overline{r}}}
\def\sq{{\overline{s}}}
\def\tq{{\overline{t}}}

\def\xq{{\overline{x}}}
\def\yq{{\overline{y}}}
\def\zq{{\overline{z}}}
%
\newdimen\Parindent
\Parindent=\parindent

\def\textindent#1{\noindent\hskip\Parindent\llap{#1\enspace }\ignorespaces}
\def\itemitem{\par\indent\hangindent2\Parindent\textindent}


\abovedisplayskip 9.0pt plus 3.0pt minus 3.0pt
\belowdisplayskip 9.0pt plus 3.0pt minus 3.0pt
\newdimen\Grenze\Grenze2\Parindent\advance\Grenze1em
\newdimen\Breite
\newbox\DpBox
\def\NewDisplay#1
#2$${\Breite\hsize\advance\Breite-\hangindent
\setbox\DpBox=\hbox{\hskip2\Parindent$\displaystyle{%
\ifx0#1\relax\else\eqno{#1}\fi#2}$}%
\ifnum\predisplaysize<\Grenze\abovedisplayskip\abovedisplayshortskip
\belowdisplayskip\belowdisplayshortskip\fi
\global\futurelet\nexttok\WEITER}
\def\WEITER{\ifx\nexttok\qed\expandafter\leftQEDdisplay
\else\leftdisplay\fi}
\def\leftdisplay{\hskip-\hangindent\leftline{\box\DpBox}$$}
\def\leftQEDdisplay{\hskip-\hangindent
\line{\copy\DpBox\hfill\lower\dp\DpBox\copy\QEDbox}%
\belowdisplayskip0pt$$\bigskip\let\nexttok=}
\everydisplay{\NewDisplay}
\newcount\GNo\GNo=0
\newcount\maxEqNo\maxEqNo=0
\def\eqno#1{%
\global\advance\GNo1
\edef\FTEST{$\fam0(\Aname.\number\GNo)$}
\ifx?#1?\relax\else
\ifnum#1>\maxEqNo\global\maxEqNo=#1\fi%
\expandafter\ifx\csname E#1\endcsname\FTEST\relax\else
\immediate\write16{E#1 hat sich geaendert!}\fi
\expandwrite\AUX{\neverexpand\ref{E#1}{\FTEST}}\fi
\llap{\hbox to 40pt{\ifx?#1?\relax\else\marginnote{E#1}\fi\FTEST\hfill}}}

\catcode`@=11
\def\eqalignno#1{\null\!\!\vcenter{\openup\jot\m@th\ialign{\eqno{##}\hfil
&\strut\hfil$\displaystyle{##}$&$\displaystyle{{}##}$\hfil\crcr#1\crcr}}\,}
\catcode`@=12

\newbox\QEDbox
\newbox\nichts\setbox\nichts=\vbox{}\wd\nichts=2mm\ht\nichts=2mm
\setbox\QEDbox=\hbox{\vrule\vbox{\hrule\copy\nichts\hrule}\vrule}
\def\qed{\leavevmode\unskip\hfil\null\nobreak\hfill\copy\QEDbox\medbreak}
\newdimen\HIindent
\newbox\HIbox
\def\setHI#1{\setbox\HIbox=\hbox{#1}\HIindent=\wd\HIbox}
\def\HI#1{\par\hangindent\HIindent\hangafter=0\noindent\leavevmode
\llap{\hbox to\HIindent{#1\hfil}}\ignorespaces}

\newdimen\maxSpalbr
\newdimen\altSpalbr
\newcount\Zaehler


\newif\ifxxx

{\catcode`/=\active

\gdef\beginrefs{%
\xxxfalse
\catcode`/=\active
\def/{\string/\ifxxx\hskip0pt\fi}
\def\TText##1{{\xxxtrue\tt##1}}
\expandafter\ifx\csname Spaltenbreite\endcsname\relax
\def\Spaltenbreite{1cm}\immediate\write16{Spaltenbreite undefiniert!}\fi
\expandafter\altSpalbr\Spaltenbreite
\maxSpalbr0pt
\gdef\alt{}
\def\\##1\relax{%
\gdef\neu{##1}\ifx\alt\neu\global\advance\Zaehler1\else
\xdef\alt{\neu}\global\Zaehler=1\fi\xdef\SigText{##1\the\Zaehler}}
\def\L|Abk:##1|Sig:##2|Au:##3|Tit:##4|Zs:##5|Bd:##6|S:##7|J:##8|xxx:##9||{%
\def\SigText{##2}\global\setbox0=\hbox{##2\relax}
\edef\TEST{[\SigText]}
\expandafter\ifx\csname##1\endcsname\TEST\relax\else
\immediate\write16{##1 hat sich geaendert!}\fi
\expandwrite\AUX{\neverexpand\ref{##1}{\TEST}}
\setHI{[\SigText]\ }
\ifnum\HIindent>\maxSpalbr\maxSpalbr\HIindent\fi
\ifnum\HIindent<\altSpalbr\HIindent\altSpalbr\fi
\HI{\marginnote{##1}[\SigText]}
\ifx-##3\relax\else{##3}: \fi
\ifx-##4\relax\else{##4}{\sfcode`.=3000.} \fi
\ifx-##5\relax\else{\it ##5\/}\fi
\ifx-##6\relax\else\ {\bf ##6}\fi
\ifx-##8\relax\else\ ({##8})\fi
\ifx-##7\relax\else, {##7}\fi
\ifx-##9\relax\else, \TText{##9}\fi\Par}
\def\B|Abk:##1|Sig:##2|Au:##3|Tit:##4|Reihe:##5|Verlag:##6|Ort:##7|J:##8|xxx:##9||{%
\def\SigText{##2}\global\setbox0=\hbox{##2\relax}
\edef\TEST{[\SigText]}
\expandafter\ifx\csname##1\endcsname\TEST\relax\else
\immediate\write16{##1 hat sich geaendert!}\fi
\expandwrite\AUX{\neverexpand\ref{##1}{\TEST}}
\setHI{[\SigText]\ }
\ifnum\HIindent>\maxSpalbr\maxSpalbr\HIindent\fi
\ifnum\HIindent<\altSpalbr\HIindent\altSpalbr\fi
\HI{\marginnote{##1}[\SigText]}
\ifx-##3\relax\else{##3}: \fi
\ifx-##4\relax\else{##4}{\sfcode`.=3000.} \fi
\ifx-##5\relax\else{(##5)} \fi
\ifx-##7\relax\else{##7:} \fi
\ifx-##6\relax\else{##6}\fi
\ifx-##8\relax\else{ ##8}\fi
\ifx-##9\relax\else, \TText{##9}\fi\Par}
\def\Pr|Abk:##1|Sig:##2|Au:##3|Artikel:##4|Titel:##5|Hgr:##6|Reihe:{%
\def\SigText{##2}\global\setbox0=\hbox{##2\relax}
\edef\TEST{[\SigText]}
\expandafter\ifx\csname##1\endcsname\TEST\relax\else
\immediate\write16{##1 hat sich geaendert!}\fi
\expandwrite\AUX{\neverexpand\ref{##1}{\TEST}}
\setHI{[\SigText]\ }
\ifnum\HIindent>\maxSpalbr\maxSpalbr\HIindent\fi
\ifnum\HIindent<\altSpalbr\HIindent\altSpalbr\fi
\HI{\marginnote{##1}[\SigText]}
\ifx-##3\relax\else{##3}: \fi
\ifx-##4\relax\else{##4}{\sfcode`.=3000.} \fi
\ifx-##5\relax\else{In: \it ##5}. \fi
\ifx-##6\relax\else{(##6)} \fi\PrII}
\def\PrII##1|Bd:##2|Verlag:##3|Ort:##4|S:##5|J:##6|xxx:##7||{%
\ifx-##1\relax\else{##1} \fi
\ifx-##2\relax\else{\bf ##2}, \fi
\ifx-##4\relax\else{##4:} \fi
\ifx-##3\relax\else{##3} \fi
\ifx-##6\relax\else{##6}\fi
\ifx-##5\relax\else{, ##5}\fi
\ifx-##7\relax\else, \TText{##7}\fi\Par}
\bgroup
\baselineskip12pt
\parskip2.5pt plus 1pt
\hyphenation{Hei-del-berg Sprin-ger}
\sfcode`.=1000
\References
}}

\def\endrefs{%
\expandwrite\AUX{\neverexpand\ref{Spaltenbreite}{\the\maxSpalbr}}
\ifnum\maxSpalbr=\altSpalbr\relax\else
\immediate\write16{Spaltenbreite hat sich geaendert!}\fi
\egroup\write16{Letzte Gleichung: E\the\maxEqNo}
\write16{Letzte Aufzaehlung: I\the\maxItemcount}}


\def\message#1{\relax}
\input xy
\xyrequire{curve}
\xyrequire{cmtip}
\xyrequire{matrix}
\xyrequire{arrow}

\def\into{\mathrel{\kern-3pt\xymatrix@=10pt{\ar@{>->}[r]&}\kern-5pt}}
\def\Auf#1{\mathrel{\kern-3pt\xymatrix@=10pt{\ar@{>>}[r]^{#1}&}\kern-5pt}}
\newdir{ >}{{}*!/-5pt/@{>}}
\SelectTips{cm}{10}
\def\cxymatrix#1{\vcenter{\xymatrix@=15pt{#1}}}
\def\inj{\ar@{ >->}}
\def\sur{\ar@{>>}}
\def\mathcite#1{\expandafter\ifx\csname#1\endcsname\relax\relax\else
\edef\next{\cite{#1}}\expandafter\dollarstrip\next\fi}
\def\dollarstrip$#1${#1}
\def\doppelpfeil{\mathop{\lower3pt\vbox{\baselineskip0pt\hbox{$\rightarrow$}\hbox{$\rightarrow$}}}}

\def\lrarrow{\mathop{\lower3pt\vbox{\baselineskip0pt\hbox{$\rightarrow$}\hbox{$\leftarrow$}}}}

\def\rho{\varrho}

\def\*{{\bf1}}
\def\0{{\bf 0}}
\fontdef BBB:bbold10,bbold7,bbold5.
\def\1{{\BBB1}}
\def\FS{{\Ss Set}}
\def\FV{{\Ss Mod}}
\def\Rel{{\Ss Rel}}

\def\Bh{{\widehat B}}

\def\cTq{{\overline\cT}}
\def\io{\leer}

\def\hat{\widehat}

\def\ITEM#1{\hangindent\Parindent\noindent{#1}}

\def\dddots{\mathinner{\mkern1mu\raise1pt\vbox{\kern7pt\hbox{.}}\mkern2mu
    \raise4pt\hbox{.}\mkern2mu\raise7pt\hbox{.}\mkern1mu}}
\def\hwidth#1{\hbox to 0pt{\hss\vrule height0pt depth 5pt width0pt$#1$\hss}}

\newcount\excount
\def\Bsp.{\noindent\advance\excount1\relax{\bf\the\excount.}\ \ignorespaces}
\let\Beispiele=\Examples
\def\Examples{\excount=0\Beispiele}


\title{Tensor envelopes of regular categories}

{\baselineskip12pt \noindent {\bf Abstract:} We extend the
calculus of relations to embed a regular category $\cA$ into a
family of pseudo\_abelian tensor categories $\cT(\cA,\delta)$
depending on a degree function $\delta$. Assume that all objects
have only finitely many subobjects. Then our results are as follows:

\item{1.} Let $\cN$ be the maximal proper tensor ideal of
$\cT(\cA,\delta)$. We show that $\cT(\cA,\delta)/\cN$ is semisimple
provided that $\cA$ is exact and Mal'cev. Thereby, we produce many new
semisimple, hence abelian, tensor categories.

\item{2.} Using lattice theory, we give a simple numerical criterion
for the vanishing of $\cN$.

\item{3.} We determine all degree functions for which
$\cT(\cA,\delta)/\cN$ is Tannakian. As a result, we are able to
interpolate the representation categories of many series of profinite
groups such as the symmetric groups $S_n$, the hyperoctahedral groups
$S_n\semidir\ZZ_2^n$, or the general linear groups $GL(n,\FF_q)$ over
a fixed finite field.

\noindent This paper generalizes work of Deligne, who first constructed the
interpolating category for the symmetric groups $S_n$.

\bigskip

\noindent
Keywords: tensor categories, semisimple categories, regular
categories, Mal'cev categories, Tannakian categories, M\"obius
function, lattices, profinite groups

\bigskip

\catcode`@=11
\def\toc#1#2{\line{\hskip20pt#1\leaders\hbox{.}\hfil#2\hskip20pt}}
\vskip0pt plus 1fil

\newread\TOC\immediate\openin\TOC=\jobname.toc
\ifeof\TOC\immediate\write16{\jobname.toc gibt es nicht!}\else
\input \jobname.toc
\fi\immediate\closein\TOC

}
\newwrite\TOC
\immediate\openout\TOC=\jobname.toc
\catcode`@=12
\vskip0pt plus 1fil
\eject

\beginsection introduction. Introduction

A category $\cA$ is called {\it regular} if it has all finite limits,
has images, and where pull\_backs preserve images. These are exactly
the prerequisites for the {\it calculus of relations}. Recall that a
relation (a.k.a. correspondence) between two objects $x$ and $y$ is a
subobject $r$ of $x\times y$. Let $s\into y\times z$ be a second
relation. Then the product $s\circ r$ of $r$ and $s$ is, by definition,
the image of $r\times_ys$ in $x\times z$. The category of relations
$\Rel(\cA)$ has the same objects as $\cA$ but with relations as
morphisms and the product of relations as composition.

In some applications, this procedure is too simplistic. For example, it
does not conform to common practice in algebraic
geometry\footnote{This example is for motivation only. Our
construction does not generalize cycle multiplication.}: let $X$, $Y$
and $Z$ be smooth complex projective varieties. Then the product of
two cycles $C\subseteq X\times Y$ and $D\subseteq Y\times Z$ is not
just the image $E$ of $C\times_YD$ in $X\times Z$. It is rather a
multiple of it (at least if $C\times_YD$ is irreducible), the factor
being the degree of the surjective morphism $e:C\times_YD\auf E$.

Guided by this example, we modify the relational product as
follows. Fix a commutative field $K$ and a map $\delta$ which assigns
to any surjective morphism $e$ of $\cA$ an element $\delta(e)$ of $K$
(its ``degree''). We define the product of $r\into x\times y$ and
$s\into y\times z$ as
$$86
sr:=\delta(e)s\circ r.
$$
where $e$ is the surjective morphism $r\times_ys\auf s\circ r$. Now,
we define a new category $\cT^0(\cA,\delta)$ as follows: it has the
same objects as $\cA$, the morphisms are formal $K$\_linear
combinations of relations, and the composition is given (on a basis)
by \cite{E86}. Of course, the degree function $\delta$ has to satisfy
certain requirements for this to work. See \cite{degreedef} for details.

The category $\cT^0(\cA,\delta)$ is only of auxiliary nature. Since
it is $K$\_linear we can enlarge it by formally adjoining direct sums
and images of idempotents (the pseudo\_abelian closure). The result is
our actual object of interest, the category $\cT(\cA,\delta)$.

This category contains in the usual way $\cA$ as a subcategory. But it
has more structure: the direct product on $\cA$ is converted into a
tensor functor on $\cT(\cA,\delta)$. It is not difficult to see that
this way, $\cT(\cA,\delta)$ is a rigid, symmetric, monoidal category (a
tensor category, for short). Loosely speaking, this means that the
tensor product has a unit element, is associative and commutative, and
that every object has a dual.

In the rest of the paper we investigate the structure of
$\cT(\cA,\delta)$. Every $K$\_linear tensor category has a maximal
proper ideal (i.e., a certain class of morphisms) which is compatible
with the tensor structure: the tensor radical $\cN$. The quotient
$\cTq(\cA,\delta)=\cT(\cA,\delta)/\cN$ is again a $K$\_linear,
pseudo\_abelian tensor category. Since its tensor radical vanishes,
$\cTq(\cA,\delta)$ has a chance to be a semisimple tensor category,
i.e., one where every object is a direct sum of simple objects. This
would entail, in particular, that $\cTq(\cA,\delta)$ is an abelian
tensor category. In our first main theorem (\cite{maintheorem}), we
show that $\cTq(\cA,\delta)$ is indeed semisimple for a large class of
categories. Moreover, we are able to determine all simple
objects. This way, we get a large number of new semisimple tensor
categories. They are non\_standard in the sense that they are not the
representation category of a (pro\_)reductive group. This was one of
the main motivations of this paper.

The precise conditions for semisimplicity are that $\cA$ is {\it
subobject finite}, {\it exact} and {\it Mal'cev.} Here, ``subobject
finite'' means that every object has only finitely many
subobjects. This is required to make all morphism spaces finite
dimensional. A regular category is exact if every equivalence relation
has a quotient while Mal'cev essentially means that all relations are
pull\_backs. These last two conditions are quite technical and it is
not clear whether they are required.

Nevertheless, the class of subobject finite, exact Mal'cev categories
has many interesting examples: the categories of finite groups, finite
rings (with or without unit), finite modules over a finite ring or,
more generally, any subobject finite abelian category, or any finite
algebraic structure containing a group operation. A particular
interesting example is the category opposite to the category of finite
sets. In that case, the construction of $\cT(\cA,\delta)$ is due to
Deligne \cite{De}.

The construction of $\cTq(\cA,\delta)$ is quite implicit since it
involves the (unknown) tensor radical $\cN$. Therefore, it is a
natural question when in fact $\cTq(\cA,\delta)$ is equal to
$\cT(\cA,\delta)$, i.e., when $\cN$ vanishes. We call degree functions
with this property {\it non\_singular}. Our second main result is a precise
numerical criterion for non\_singularity. The only assumption on $\cA$
is subobject finiteness. For a surjective $\cA$\_morphism
$e:x\auf y$ we define the number
$$
\omega_e:=\sum_{w\subseteq x\atop e(w)=y}\mu(w,x)\delta (w\auf y)\in K
$$
where $\mu$ is the M\"obius function on the lattice of subobjects
of $x$. A surjective morphism $e$ is {\it indecomposable} if $e$ is
not an isomorphism and if any factorization $e=e'e''$ into surjective
morphisms implies that one of $e'$ or $e''$ is an isomorphism.  The
criterion is that {\it $\delta$ is non\_singular if and only if
$\omega_e\ne0$ for all indecomposable $e$.}

With this criterion it is very easy to compute the singular degree
functions in many cases. For example, the degree functions of the
category $\cA=\FS^{\rm op}$ are parametrized by one number $t\in
K$. The corresponding degree function is singular precisely when
$t\in\NN$, recovering a result of Deligne. Similarly, for
$\cA=\FV_{\FF_q}$, the category of finite $\FF_q$\_vector spaces, the
singular parameters are precisely the powers $q^n$ with $n\in\NN$. On
the more abstract side, we can show that there always {\it exists} a
non\_singular degree function provided that $\cA$ is exact and
protomodular. The latter condition on $\cA$ is stronger than Mal'cev
but holds for all the examples mentioned above.

The best known semisimple tensor categories are the representation
categories of pro\_reductive groups (so called {\it Tannakian
categories}\footnote{At least if $K$ is algebraically closed. Assume
this from now on.}). Thus it is a natural problem to determine degree
functions $\delta$ for which $\cTq(\cA,\delta)$ is Tannakian.  Our
third main result answers this question roughly as follows: assume
$\cA$ is a subobject finite, regular category and that $K$ is
algebraically closed of characteristic zero. Then $\cTq(\cA,\delta)$
is Tannakian if and only if $\delta$ is adapted to a uniform functor
$P:\cA\pfeil\FS$. In this case, $\cTq(\cA,\delta)\cong\|Rep|(G,K)$
where $G$ is the profinite group of automorphisms of $P$ (see
Definitions~\ncite{unidef} and~\ncite{adapteddef} concerning
``uniform'' and ``adapted'').

We don't know of a construction of uniform functors in general but, in
examples, it is not difficult to come up with many of them. More
precisely, for certain categories $\cA$ we are able to construct
sufficiently many uniform functors $P_i$, $i\in I$, such that the
corresponding adapted degree function $\delta_i$ are Zariski\_dense in
the space of all degree functions. Let $G_i:=\|Aut|(P_i)$ be the
associated group. Since $\|Rep|(G_i,K)$ is a quotient of
$\cT(\cA,\delta_i)$ we say that {\it $\cT(\cA,\delta)$ interpolates
the categories $\|Rep|(G_i,K)$, $i\in I$}.

Let for example $\cA=\FS^{\rm op}$. As already mentioned, it has a
one\_parameter family of degree functions $\delta_t$. It turns out
that $\cTq(\cA,\delta_t)\equiv\|Rep|(S_n,K)$ when $t=n\in\NN$
(coincidentally(?) precisely the parameters for which $\delta_t$ is
singular). Thus $\cT(\cA,\delta_t)$ interpolates the representation
categories of the symmetric groups $S_n$, $n\in\NN$ (that was
Deligne's starting point). Similarly, we find a category
$\cT(\cA,\delta_t)$ which interpolates the representation categories
of $GL(n,\FF_q)$, $n\in\NN$, $q$ fixed. Other examples include the
family of wreath products $S_n\wr G$, for $G$ a fixed finite group, or
even the infinite wreath product $S_{n_1}\wr S_{n_2}\wr S_{n_3}\ldots$
and many more. We hope that our construction gives rise to a
simultaneous treatment of the representations of the $G_i$, in the
same way as the representations of the symmetric groups are best
studied simultaneously.

The paper concludes with two appendices. In the first one, we give a
very brief introduction to protomodular and Mal'cev categories. As
already mentioned, we need ``Mal'cev'' for proving semisimplicity and
``protomodular'' for the existence of a non\_singular degree
function. In the second appendix, we use the Mal'cev property to
compute degree functions.

We have tried to enhance our theory by including a fair number of
examples. In addition to some isolated ones, the paper contains five
more extensive blocks of examples. They cover regular categories
(section~\cite{Regular}), degree functions (section \cite{CoCor}),
singular degree functions (section \cite{determinant}), interpolation
of Tannakian categories (section \cite{Tannakian}), and
protomodular/Mal'cev categories (Appendix A).

The present work owes its existence to the paper \cite{De} of Deligne
where he constructs $\cT(\cA,\delta)$ in the case $\cA=\FS^{\rm op}$.
His construction is carried out using different building blocks but
the result is the same. Also the backbone of the proofs of our three
main results is taken from Deligne's paper. We just added some more
flesh to it. The main novelty of the present paper is probably the
identification of the Mal'cev condition as being the key for the
semisimplicity proof and the numerical non\_singularity criterion in
terms of M\"obius functions.

Finally, it should be mentioned that this paper has a predecessor,
\cite{KnTen}, where the theory is is carried out in the special case
of abelian categories. One of my motivations for the present paper was
to bring Deligne's case $\FS^{\rm op}$ and the case of abelian
categories under a common roof.

\beginsection Regular. Regular categories

Regular categories have been introduced by Barr, \cite{Barr}, but the
extent limits are supposed to exist in their definition varies from
author to author. In this section we make our notion of regularity
precise and set up some other terminology.

Let $\cA$ be a category. Monomorphisms in $\cA$ will henceforth be
called {\it injective} and will be indicated by the arrow
``$\into$''. Two injective morphisms $f:u\into x$ and $f':u'\into x$
with the same target are {\it equivalent}, $f\approx f'$, if there
exists an isomorphism $g:u\pf\sim v$ with $f=f'g$. A {\it subobject}
of $x$ is an equivalence class of injective morphisms. We denote the
class of subobjects of $x$ by $\|sub|(x)$. In most of this paper, we
are going to assume that $\|sub|(x)$ is a set ({\it well\_powered}) or
even finite ({\it subobject finite}) for all $x$. The set $\|sub|(x)$
has the structure of a poset: in the notation above, we say $f\le f'$
(or just $u\subseteq u'$) if there is a morphism $g$ with $f=f'g$. The
morphism $g$ is injective and unique. Hence $f\le f'$ and $f'\le f$
imply $f\approx f'$.

The image, $\|image|(f)$, of any morphism $f:x\pfeil y$ is the
(absolutely) smallest subobject of $y$ through which $f$
factorizes. Clearly, the image may or may not exist. The morphism $f$
will be called {\it surjective} (or, more traditionally, an {\it
extremal epimorphism}) if $\|image|(f)=y$. A surjective morphism will
be indicated by the arrow ``$\auf$''.

\resetitem

\Definition. A category $\cA$ is {\it complete and regular} if
\item{\bf R0} $\cA$ is well\_powered, i.e., $\|sub|(x)$ is a set for every
object $x$.
\item{\bf R1}$\cA$ has all finite limits. In particular, it has a terminal
object denoted by $\*$.
\item{\bf R2}Every morphism has an image.
\item{\bf R3}The pull\_back of a surjective morphism along any morphism
is surjective.\Par

\Remarks: 1. The first axiom, {\bf R0}, is non\_standard and is only
thrown in for convenience.

\noindent 2. Axiom~{\bf R2} (together with {\bf R1}) implies that
every morphism can be factorized as $f=me$ where $m$ is injective and
$e$ is surjective. This factorization is essentially unique. Moreover,
the classes of surjective and injective morphisms are closed under
composition and their intersection consists of the isomorphisms.

\noindent 3. Usually, only regular epimorphisms are called
surjective. In particular, {\bf R3} is stated only for regular
epimorphisms. One can show (see e.g. \cite{Borceux}~\S2.2) that, in
the presence of {\bf R1}--{\bf R3}, the concepts ``extremal
epimorphism'', ``strong epimorphism'', and ``regular epimorphism'' are
all the same.

\medskip

For any category $\cA$ let $\cA^\io$ be the category obtained by
formally adjoining a (new) absolutely initial object $\io$. More
precisely, an object of $\cA^\io$ is either an object of $\cA$ or
equal to $\io$. The morphisms between objects of $\cA$ stay the same,
for every object $x$ of $\cA^\io$ there a unique morphism from $\io$
to $x$ and no morphism from $x$ to $\io$ unless $x=\io$.

\Definition. A non\_empty category $\cA$ is {\it regular} if $\cA^\io$
is complete and regular.

\noindent In down to earth terms, this means:

\Proposition partreg. A category $\cA$ is regular if it satisfies {\bf
R0}, {\bf R2}, {\bf R3} above and if {\bf R1} is replaced by:
\itemitem{\bf R1.1}$\cA$ has a terminal object $\*$.
\itemitem{\bf R1.2}For every commutative diagram
$$49
\cxymatrix{y\ar[r]\ar[d]&v\ar[d]\\u\ar[r]&x}
$$
the pull\_back $u\times_xv$ exists.
\itemitem{\bf R1.3}The pull\_back of a surjective morphism by an arbitrary
morphism exists.\Par

\Proof: First observe that the inclusion of $\cA$ in $\cA^\io$
preserves and reflects limits, injective morphism, and equality of
subobjects.The same holds then for images and surjective
morphisms. Then one checks easily:
$$0
\eqalignno{
&{\bf R0}^\io&\Leftrightarrow{\bf R0}\cr
&{\bf R1}^\io&\Leftrightarrow{\bf R1.1}\ {\rm and}\ {\bf R1.2}\cr
&{\bf R2}^\io&\Leftrightarrow{\bf R2}\cr
&{\bf R3}^\io&\Leftrightarrow{\bf R1.3}\ {\rm and}\ {\bf R3}\cr
}
$$
where ${\bf R}i^\io$ means axiom ${\bf R}i$ for $\cA^\io$.\qed

\excount=0
\Remarks: \Bsp. Recall that a {\it cone} of a diagram $D:\cD\pfeil\cA$
is an object $x$ together with morphisms $f_d:x\pfeil D(d)$ which
satisfy the obvious commutation relations. Call a diagram {\it
bounded} if it has a cone. Then {\bf R1.1} and {\bf R1.2} are
equivalent to the following completeness statement: {\it every bounded
finite diagram has a limit.} This implies in particular that every
regular category with an initial object is complete.

\Bsp. Many authors define regular categories to be complete. We opted
for our present terminology mainly for two reasons. First, it
accommodates some (for me) important examples, namely the category of
(non\_empty) affine spaces and the categry of free actions of a
group. Socondly, it has also conceptual advantages. See, e.g., the
decomposition \cite{iddecomp} below.

\medskip

\noindent In the following we use freely the embedding of $\cA$ into
$\cA^\io$ in the way that $\io$ stands for all non\_existent limits.

\Examples: Regular categories, even complete ones are abundant. The
category of models of any equational theory is complete and
regular. This includes the categories of sets, lattices, groups,
rings, etc. The category of compact Hausdorff spaces is complete
regular as is every abelian category. Also the category {\it opposite}
to the category of sets is regular.

One reason for the abundance is that the concept of regular
categories enjoys many permanence properties. The list below is not
exhaustive. In the following let $\cA$ be a regular category.

\Bsp. Let $\cB$ be a full subcategory of $\cA$ which contains the
   terminal object and is closed under products and subobjects. Then
   $\cB$ is regular. This applies in particular to the category of
   finite models of an equational theory: finite sets, finite
   lattices, finite groups, etc.

\Bsp. Let $\cD$ be a small category. Then the category of all
   functors $\cD^{\rm op}\pfeil\cA$ is regular (a {\it diagram
   category}). Examples are the categories of all arrows $x\pfeil y$
   in $\cA$ or the category of all objects equipped with a $G$\_action
   where $G$ is a fixed group.

\Bsp. Fix an object $s$ of $\cA$. Then the category $\cA/s$ of all
   ``$s$\_objects'', i.e., all arrows $x\pfeil s$, a so\_called {\it
   slice category}, is regular. This is one of the main mechanisms to
   obtain regular categories which are not pointed, i.e., do not
   possess a zero object.

\Bsp. For a fixed object $s$ the category $\cA\mod s$ of all {\it
   surjective} morphisms $x\auf s$ is regular. A prime example is the
   category of (non\_empty) affine spaces for a fixed field $k$. This
   category is equivalent to the category of finite dimensional
   $k$\_vector spaces equipped with a {\it non\_zero} linear form (the
   equivalence is given by taking the dual of the space of affine
   functions). As opposed to the previous constructions the categories
   produced by this one are usually not complete even if $\cA$ is. In
   fact, this is is one of our main motivations for our notion of
   regularity.

\Bsp. Again, fix an object $p$ of $\cA$. Then the category
   $p\backslash\cA$ of all arrows $p\pfeil x$ (a {\it coslice
   category}) is regular. If $p=\*$ then $\cB$ is the category of all
   pointed objects.

\Bsp. Fix an object $p$ and consider the full subcategory
   $p\pfeil\cA$ of objects $x$ for which there {\it
   exists} a morphism $p\pfeil x$. It is regular if $p$ is projective
   in $p\pfeil\cA$, i.e., if for any surjective
   morphism $u\auf v$ such that there is a morphism $p\pfeil u$ every
   morphism $p\pfeil v$ can be lifted to $p\pfeil u$. Take, e.g., for
   $\cA$ the category {\it opposite} to the category of $G$\_sets ($G$
   a fixed group). Let $p=G$ with left regular action. Then
   $p\pfeil\cA$ is the opposite category of the
   category of $G$\_sets with {\it free} $G$\_action. Again, this
   example is only regular and not complete.

\Bsp. A combination of slice and coslice category is the category of
   $s$\_points $\|Pt|_s\cA=s\backslash(\cA/s)=(s\backslash\cA)/s$. It
   is the category of all triples $(x,e,d)$ where $e:x\pfeil s$ is a
   morphism and $d:s\pfeil x$ is a section of $e$. Its main virtue is
   that it is a pointed category.

\medskip

In every regular category there is also the notion of a {\it quotient
object of $x$}. It is an equivalence class of surjective morphisms
with domain $x$. The {\it kernel pair} of an quotient object $x\auf y$
is the double arrow $x\times_yx\doppelpfeil x$. It determines the
quotient object uniquely since $y$ is the coequalizer of its kernel
pair. In other words, every quotient object is encoded by the
subobject $x\times_yx$ of $x\times x$.

The kernel pair is an example of an equivalence relation. In general,
an {\it equivalence relation on $x$} is a subobject $r$ of $x\times x$
which is reflexive (i.e., contains the diagonal), symmetric (i.e., is
invariant under exchanging the two factors of $x\times x$, and
transitive (i.e., the morphism $r\times_xr\pfeil x\times x$ factorizes
through $r$. In general, not every equivalence relation is the kernel
pair $x\times_yx\into x\times x$ of a quotient object. If it is then
it is called {\it effective}. Thus, there is a bijection between
quotient objects of $x$ and effective equivalence relations on $x$.

\Definition Exactdef. A category is {\it exact} if it is regular and
if every equivalence relation is effective.

An object $y$ is a {\it subquotient} of an object $x$ if $y$ is a
quotient of a subobject of $x$, i.e., if there is a diagram
$\xymatrix@=10pt{x&u\ar@{ >->}[l]\ar@{>>}[r]&y}$. In that case, we
write $y\preceq x$. If we can find such a diagram such that at least
one of the two arrows is not an isomorphism then this is denoted by
$y\prec x$. In that case we call $y$ a {\it proper} subquotient of
$x$.

\Lemma. The relations ``$\preceq$'' and ``$\prec$'' are transitive.

\Proof: Assume $z\preceq y\preceq x$. Then, we get a diagram
$$52
\cxymatrix{
x''\inj[r]\sur[d]&x'\inj[r]\sur[d]&x\\
y'\inj[r]\sur[d]&y\\
z
}
$$
where the square is a pull\_back showing $z\preceq x$. If
both $x''\pfeil x$ and $x''\pfeil z$ were isomorphisms then all
morphisms in diagram \cite{E52} were isomorphism showing that
``$\prec$'' is transitive, as well.\qed

\Lemma partialorder0. Let $x$ be an object such that $\|sub|(x)$
satisfies the descending chain condition and $\|sub|(x\times x)$
satisfies the ascending chain condition. Then there is no infinite
chain
$$31
x\succ x_1\succ x_2\succ x_3\succ\ldots
$$
In particular, we have $x\not\prec x$.

\Proof: Since the quotient object $x\auf y$ is determined by the
subobject $x\times_yx\into x\times x$ the ascending chain condition
for $\|sub|(x\times x)$ implies the descending chain condition for
quotients of $x$. Let $z\into x$ be any subobject. Since
$\|sub|(z\times z)\subseteq\|sub|(x\times x)$ we see that every
subobject of $x$ satisfies the descending chain condition on quotient
objects.

The chain \cite{E31} gives rise to the diagram
$$
\cxymatrix{\cdots&z_{14}\inj[r]\sur[d]&z_{13}\inj[r]\sur[d]&
z_{12}\inj[r]\sur[d]&x_1\\
\cdots&z_{24}\inj[r]\sur[d]&z_{23}\inj[r]\sur[d]&x_2\\
\cdots&z_{34}\inj[r]\sur[d]&x_3\\
\cdots&x_4\\
\dddots
}
$$
where all squares are pull\_backs. By the descending chain condition
for $\|sub|(x)$ there is a bound $N>0$ such that $z_{1\,n+1}\pf\sim
z_{1\,n}$ for all $n\ge N$. Then all morphisms in the second row
$z_{2\,n+1}\into z_{2n}$ are also surjective, hence isomorphisms, for
$n\ge N$. We conclude that all horizontal arrows $z_{i\,n+1}\into
z_{in}$ are isomorphisms for $n\ge N$ and $1\le i\le n$. But then
we get an infinite chain of quotients $z_{1N}\auf x_N\auf
x_{N+1}\auf\ldots$ showing that $x_{n+1}$ can't be a proper
subquotient of $x_n$ for $n>\!\!>0$.\qed

\beginsection CoCor. The construction of the tensor envelope $\cT(\cA,\delta)$

First, we recall the classical calculus of relations. Let $\cA$ be a
complete regular category. A {\it relation} (a.k.a. {\it
correspondence}) between $x$ and $y$ is by definition a subobject $r$
of $x\times y$. Let $s$ be a relation between $y$ and $z$. Then the
{\it product} of $r$ and $s$ is defined as
$$
s\circ r:=\|image|(r\times_ys\pfeil x\times z)
$$
or, as a diagram,
$$37
\cxymatrix{
&&\hwidth{r\times_ys}\ar[dl]\ar[dr]\ar@{.>>}[d]^e\\
&r\ar[dl]\ar[dr]&\hwidth{s\circ r}\ar@{.>}[dll]\ar@{.>}[drr]&s\ar[dl]\ar[dr]\\
x&&y&&z
}
$$
The regularity of $\cA$ (or, more precisely, axiom {\bf R3}) ensures
that this product is associative. This way, one can define a new
category $\Rel(\cA)$ with the same objects as $\cA$ but with relations
as morphisms.

The construction of $\Rel(\cA)$ completely ignores the structure of
the surjective morphism $e$ in diagram \cite{E37}. Our main
construction can be roughly described as replacing $e$ be a numerical
factor, its ``degree'' or ``multiplicity''. To carry this out we need
to consider linear combinations of relations which actually enlarges
the scope of the construction: since there is now a zero morphism not
all pull\_backs have to exist. They are just set to zero. Here are the
precise definitions:

\Definition degreedef. Let $\fE(\cA)$ be its class of surjective
morphisms of a (just) regular category $\cA$ and let $K$ be a
commutative ring. Then a map $\delta:\fE(\cA)\pfeil K$ is called a
{\it degree function} if
\medskip
\item{\bf D1} $\delta(1_x)=1$ for all $x$.
\item{\bf D2} $\delta(\eq)=\delta(e)$ whenever $\eq$ is a pull\_back
of $e$.
\item{\bf D3} $\delta(e\,\eq)=\delta(e)\,\delta(\eq)$ whenever $e$ can
be composed with $\eq$.

\Examples: The degree functions in the following examples can be
determined by simple ad\_hoc arguments. Observe however that in
Appendix~B we have proved some general statements on the computation
of degree functions which cover most of the examples below.

\Bsp. The morphism $\io\pfeil\io$ of $\cA^\io$ is a pull\_back of {\it
every} morphism. Thus, the only degree function on $\cA^\io$ is the
trivial one: $\delta\equiv1$. Hence, it is {\it not} possible to reduce
to the complete case by simply replacing $\cA$ with $\cA^\io$.

\Bsp. [Deligne's case] By the same reason, all degree functions on the
category of (finite) sets are trivial. On the other hand, if $\cA$ is
the category {\it opposite} to the category of finite sets then the
surjective morphisms of $\cA$ are the injective maps in $\FS$. In that
case all degree functions are of the form
$$57
\delta(e:A\into B)=t^{|B\setminus e(A)|}
$$
where $t\in K$ is arbitrary.

\Bsp. Let $\cA$ be the category of finite dimensional $k$\_vector spaces
  where $k$ is some field. Then all degree functions are of the form
$$56
\delta(e:U\auf V)=t^{\|dim|\|ker|e}
$$
where $t\in K$ is arbitrary.

\Bsp. More generally, let $\cA$ be an abelian category in which every
   object is of finite length. Let $S$ be the class of simple
   objects. For an object $x$ let $\ell_s(x)$ be the multiplicity of
   $s\in S$ in $x$. Then all degree functions are of the form
$$55
\delta(e:x\auf y)=\prod_{s\in S}t_s^{\ell_s(\|ker|e)}
$$
   where the parameters $t_s\in K$ are arbitrary.

\medskip

Now we define a $K$\_linear category as follows:

\Definition. Let $\cA$ be a regular category, $K$ a commutative ring,
and $\delta$ a $K$\_valued degree function on $\cA$. Then the
category $\cT^0(\cA,\delta)$ is defined as follows:
\medskip
\ITEM{The \it objects} of $\cT^0(\cA,\delta)$ are those of $\cA$. If
an object $x$ of $\cA$ is regarded as an object of $\cT^0$ then we will
denote it by $[x]$.\Par
\ITEM{The \it morphisms} from $[x]$ to $[y]$ are the formal
$K$\_linear combinations of relations between $x$ and $y$. If $x\times
y$ does not exist then $\|Hom|_{\cT^0}([x],[y])=0$.\Par
\ITEM{The \it composition} of $\cT^0$\_morphisms is defined on a basis
as follows: let $r\into x\times y$ and $s\into y\times z$ be
relations. Then their composition is (in the notation of \cite{E37})
$$35
sr:=\cases{\delta(e)\ s\circ r&if
$r\times_ys$ exists\cr0&otherwise\cr}.
$$
\Par

\Remark: If $\cA$ is complete regular and $\delta\equiv1$ then
$\cT^0(\cA,\delta)$ is just the $K$\_linear hull of $\Rel(\cA)$.

\medskip

To facilitate further computations, we reformulate and extend the
product formula \cite{E35}. First, we adopt the following notation: if
$x$ and $y$ are objects of $\cA$ and $f:r\pfeil x\times y$ is any
$\cA^\io$\_morphisms (i.e., $f$ may not be injective and $r$ may be
$\io$) with image $\rq$ then we define the $\cT^0$\_morphism
$\<f\>:[x]\pfeil[y]$ as
$$
\<f\>:=\cases{\delta (r\auf\rq)\ \rq&if $r\ne\io$,\cr0&if $r=\io$.}
$$

\Lemma. Let $x$, $y$, $z$ be objects of $\cA$ and $r\pfeil x\times y$
and $s\pfeil y\times z$ be $\cA^\io$\_morphisms. Then
$$36
\<s\pfeil y\times z\>\<r\pfeil x\times y\>=\<r\times_ys\pfeil x\times z\>
$$

\Proof: If one of $r$ or $s$ equals $\io$ then $r\times_ys=\io$, and
both sides of \cite{E36} are zero. So assume $r,s\ne\io$. Let $\rq$
be the image of $r\pfeil x\times y$ and $\sq$ the image of $s\pfeil
y\times z$. Then we obtain the following diagram (in $\cA^\io$)
$$
\xymatrix@R=10pt@C=13pt{
&&t\ar[ld]\ar[rd]\\
&\tilde r\ar[ld]\ar[rd]&&\tilde s\ar[ld]\ar[rd]\\
r\ar@{>>}[rd]&&\tq\ar[ld]\ar[rd]&&s\ar@{>>}[ld]\\
&\rq\ar[ld]\ar[rd]&&\sq\ar[ld]\ar[rd]\\
x&&y&&z\\
}
$$
where all squares are pull\_backs. Axiom {\bf R3} implies that the two
morphisms $t\pfeil\tilde s$ and $\tilde s\pfeil \tq$ are
surjective. Thus $t=\io$ implies $\tq=\io$, in which case both sides of
\cite{E36} are zero. So, assume $t\ne\io$. Then we
get
$$
\eqalign{\<s\pfeil y\times z\>\<r\pfeil x\times y\>
&=\delta (r\pfeil\rq)\,\delta (s\pfeil\sq)\<\tq\pfeil x\times z\>
\ {\buildrel{\bf D2}\over=}\cr
&=\delta (t\pfeil\tilde s)\delta (\tilde s\pfeil\tq)
\<\tq\pfeil x\times z\>\ {\buildrel{\bf D3}\over=}\ \<t\pfeil x\times z\>.\cr}
$$\qed

Now we can prove:

\Theorem iscategory. Let $\cA$ be a regular category, $K$ a
commutative ring, and $\delta :\fE(\cA)\pfeil K$ a degree
function. Then $\cT^0(\cA,\delta)$ is a category.

\Proof: Condition {\bf D1} makes sure that the diagonal relation
$x\pfeil x\times x$ is an identity morphism in $\cT^0$. It
remains to show that composition is associative. Let $F$, $G$, and $H$
be the $\cT^0$\_morphisms corresponding to relations $r\into
x\times y$, $s\into y\times z$, and $t\into z\times u$,
respectively. Then
$$
\eqalign{(HG)F&=\<s\times\limits_z t\pfeil y\times u\>\cdot
\<r\pfeil x\times y\>\buildrel\mathcite{E36}\over=\cr
&=\Big\<r\times_y(s\times\limits_z t)\pfeil x\times u\Big\>
=\Big\<(r\times_ys)\times\limits_z t\pfeil x\times u\Big\>
\buildrel\mathcite{E36}\over=\cr
&=\<t\pfeil z\times u\>\cdot\<r\times\limits_y s\pfeil x\times z\>=
H(GF).\cr}
$$\qed

The category $\cT^0$ is only of auxiliary nature, our main interest
being its {\it pseudo\_abelian closure} $\cT(\cA,\delta)$. Recall that
a category is pseudo\_abelian (or also Karoubian) if it is additive
and every idempotent has an image. We give a brief description of how
to construct $\cT$. For details, see, e.g., \cite{AKO}~\S1.

The pseudo\_abelian closure of $\cT^0$ is constructed in two
steps. First one forms the additive closure $\cT'$ of $\cT^0$. Its
objects are formal direct sums $\oplus_{i=1}^n[x_i]$. Morphisms are
matrices of $\cT^0$\_morphisms. Observe that the empty direct sum
($n=0$) is allowed and provides a zero object.

The category $\cT(\cA,\delta)$ is now the idempotent closure of
$\cT'$: the objects of $\cT$ are pairs $(X,p)$ where $X$ is an object
of $\cT'$ and $p\in\|End|(X)$ is idempotent. The morphism space
between $(X,p)$ and $(X',p')$ is $p'\|Hom|_{\cT'}(X,X')p$. This
construction shows, in particular, that $\cT^0$ is a full subcategory
of $\cT$.

The category $\cA$ is a subcategory of $\cT^0$. In fact, for an
$\cA$\_morphism $f:x\pfeil y$ let $[f]\into x\times y$ be its
graph. Then one checks easily that $x\mapsto[x]$ and $f\mapsto[f]$
defines an embedding $\cA\pfeil\cT^0$. Since $\cT(\cA,\delta)$ has a
zero object, this embedding can be extended to a functor
$\cA^\io\pfeil\cT(\cA,\delta)$ by defining $[\io]=0$.

The direct product turns $\cA^\io$ into a symmetric monoidal
category. This induces a $K$\_linear tensor product on
$\cT(\cA,\delta)$ by defining
$$
[x]\otimes[y]:=[x\times y].
$$
This tensor product is functorial: for relations $\<r\into x\times
x'\>$, $\<s\into y\times y'\>$ let
$$53
\<r\>\otimes\<s\>:[x]\otimes[y]\pfeil[x']\otimes[y']
$$
be the relation
$$
r\times s\Pfeil(x\times x')\times(y\times
y')\Pf\sim(x\times y)\times(x'\times y').
$$
The unit object is $\1=[\*]$. We claim that the tensor product is
rigid, i.e., every object $X$ has a dual $X^\vee$. It suffices to
prove this for objects of the form $X=[x]$. But then $X$ is even
selfdual with evaluation morphism $\|ev|:[x]\otimes[x]\pfeil\1$ and
coevaluation morphism $\|ev|^\vee:\1\pfeil[x]\otimes[x]$ represented by
$$47
\cxymatrix{
&x\ar@{ >->}[ld]_{\rm diag}\ar[rd]&\\
x\times x&&{\*}
}
\qquad
\cxymatrix{
&x\ar[ld]\ar@{ >->}[rd]^{\rm diag}&\\
\*&&x\times x
}
$$

In a rigid monoidal category every morphism
$F:X\pfeil Y$ has a
{\it transpose}
$$
F^\vee:Y^\vee\Pfeil Y^\vee\otimes X\otimes X^\vee\Pf{1\otimes f\otimes
1}Y^\vee\otimes Y\otimes X^\vee\Pfeil X^\vee.
$$
Concretely, if $F:[x]\pfeil [y]$ is represented
by $r\pfeil x\times y$ then $F^\vee:[y]=[y]^\vee\pfeil
[x]^\vee=[x]$ is represented by the transposed relation $r\pfeil
x\times y\pf\sim y\times x$.

Every tensor category is linear over the endomorphism ring of
$\1$. Therefore, the following statement is evident but crucial.

\Proposition Endo1. $\|End|_{\cT(\cA,\delta)}(\1)$ is the free
$K$\_module with basis $\|sub|(\*)$. Multiplication is given given by
intersection (with the convention that $u\cdot v=0$ if $u\cap v$ does
not exist).

\noindent For any $u\in\|sub|(\*)$ let $\cA_u\subseteq\cA$ be the full
subcategory whose objects are those $x$ such that
$\|image|(x\pfeil\*)=u$. This is again a regular category with
terminal object $u$. The degree function $\delta$ on $\cA$ restricts to
a degree function $\delta_u$ on $\cA_u$.

\Theorem iddecomp. Assume that $\|sub|(\*)$ is finite. Then
$\cT(\cA,\delta)$ is tensor equivalent to the product of the
categories $\cT(\cA_u,\delta_u)$ where $u$ runs through $\|sub|(\*)$.

\Proof: Let $A:=\|End|_{\cT}(\1)$. Then \cite{Endo1} implies that $A$
has a basis of orthogonal idempotents $p_u$ such that
$$54
p_uv=\cases{p_u&if $u\le v$,\cr0&otherwise\cr}
$$
(see \S\cite{Moebius} below for details). Let $p_u\cT$ be the category
with the same objects as $\cT$ but with
$$
\|Hom|_{p_u\cT}(X,Y)=p_u\|Hom|_{\cT}(X,Y).
$$
It is easy to check that this is again a pseudo\_abelian. Moreover the
functor $X\mapsto (X)_u$, $F\mapsto (p_uF)_u$ is an equivalence of
tensor categories $\cT\Pf\sim\prod_up_u\cT$. It remains to show that
$p_u\cT$ is equivalent to $\cT(\cA_u,\delta_u)$. Since $p_u\cT$ is the
pseudo\_abelian closure of $p_u\cT^0$ is suffices to prove that
$p_u\cT^0$ is equivalent to $\cT^0(\cA_u,\delta_u)$.

We claim that we can define a functor
$\Phi:\cT^0(\cA_u,\delta_u)\pfeil p_u\cT^0$ by sending the object
$[x]$ to itself and a morphism $F$ to $p_uF$. The only problem is for
relations $r$ and $s$ in $\cA_u$ such that $s\circ r$ is not in
$\cA_u$. In that case let $v:=\|image|(s\circ r\pfeil \1)\subset
u$. Now, according to formula \cite{E53}, the action of
$u\in\|End|_\cT(\1)$ on a relation $r$ is given by
$$
u\<r\>=\<u\times r\>.
$$
Thus $p_u\<s\circ r\>=p_u\<v\times s\circ r\>=p_uv\<s\circ r\>=0$
by \cite{E54} proving the claim.

Now we show that $\Phi$ is a tensor equivalence. Let $x$ be any object
of $\cA$. We claim that $i:u\times x\into x$ induces an isomorphism
$[u\times x]\pf\sim[x]$ in $p_u\cT$. Indeed, $i^\vee i=1_{[u\times
x]}$ even in $\cT$. Conversely,
$p_u1_{[x]}=p_u\<x\>=p_uu\<x\>=p_u\<u\times x\>=p_uii^\vee$ which
proves the claim.

Put $v:=\|image|(x\pfeil\*)$. If $u\subseteq v$ then $u\times x$ is an
object of $\cA_u$. If $u\not\subseteq v$ then
$p_u1_{[x]}=p_uv1_{[x]}=0$, hence $[x]=0$ in $p_u\cT^0$. This shows
that every object of $p_u\cT^0$ is isomorphic to an object in the image
of $\Phi$.

Let now $x$ and $y$ be two objects of $\cA_u$. Then
$$
\|Hom|_{\cT}([x],[y])=\|Hom|_{\cT^0(\cA_u,\delta_u)}([x],[y])\oplus C
$$
where $C$ is spanned by all relations $r$ with
$\|image|(r\pfeil\*)\subset u$. This shows that
$$
\|Hom|_{\cT^0(\cA_u,\delta_u)}([x],[y])\Pf{p_u\cdot}\|Hom|_{p_u\cT}([x],[y])
$$
is an isomorphism, completing the proof that $\Phi$ is an equivalence
of categories. Finally it is a tensor equivalence since
$[x\times^{\cA_u}y]\cong[x\times^\cA y]$ in $p_u\cT^0$ for all objects
$x,y$ of $\cA_u$.\qed

\noindent For our purposes, the preceding theorem allows us to assume
without loss of generality that $\*$ has no proper subobject or,
equivalently, that the endomorphism ring of $\1$ is $K$. This is
one of the main reasons for our definition of regular categories.

\beginsection radical. The radical of a tensor category

In this section we review some general facts about tensor
categories. Details can be found, e.g. in \cite{AKO}. Let $\cT$ be an
arbitrary pseudo\_abelian tensor category and denote the commutative
ring $\|End|_\cT(\1)$ by $K$.

Let $\cI$ be a map which assigns to any two objects $X$ and $Y$ of
$\cT$ a subspace $\cI(X,Y)$ of $\|Hom|_\cT(X,Y)$. Then $\cI$ is called
a {\it tensor ideal} if

\item{a)} it is closed under arbitrary left and right
multiplication, i.e., for all diagrams $W\pf f X\pf g Y\pf h Z$
holds: if $g\in\cI(X,Y)$ then $hgf\in\cI(W,Z)$ and

\item{b)} it is closed under tensor products, i.e., for all morphisms
$f:X\pfeil Y$ and objects $Z$ holds: if $f\in\cI(X,Y)$ then $f\otimes
1_Z\in\cI(X\otimes Z,Y\otimes Z)$. 

\noindent Given a tensor ideal $\cI$, it is possible to define a
tensor category $\cT/\cI$. Its objects are the same as those
of $\cT$ but the morphisms are:
$$
\|Hom|_{\cT/\cI}(X,Y):=\|Hom|_\cT(X,Y)/\cI(X,Y).
$$
In fact, property a) makes sure that composition of morphisms can be
pushed down to $\cT/\cI$. Property b) does the same for morphisms
between tensor products. The category $\cT/\cI$ is clearly
additive. For pseudo\_abelian we need a further condition.

\Lemma. Assume $K$ is an Artinian ring and that all
$\|Hom|_\cT$\_spaces are finitely generated $K$\_modules. Then
$\cT/\cI$ is also pseudo\_abelian.

\Proof: Follows from the following well known fact: let $A\auf B$ be a
surjective homomorphism between Artinian rings. Then every idempotent
of $B$ can be lifted to an idempotent of $A$.\qed

Using the isomorphism
$$
\iota_{XY}:\|Hom|_\cT(\1,X^\vee\otimes Y)\Pf\sim\|Hom|_\cT(X,Y).
$$
we have (\cite{AKO}~6.1.5.)
$$40
\cI(X,Y)=\iota_{XY}(\cI(\1,X^\vee\otimes Y)).
$$
This implies, in particular, that $\cI=0$ if and only if $\cI(\1,Y)=0$
for all $Y$.

Let $f:X\pfeil X$ be an endomorphism in $\cT$. The trace, $\|tr|f$, of $f$
is the composition
$$
\1\Pf\delta X\otimes X^\vee\Pf{f\otimes1_{X^\vee}}X\otimes X^\vee\Pf\sim
X^\vee\otimes X\Pf{\rm ev}\1.
$$
The trace is an element of $K=\|End|_\cT(\1)$. Now we define
the {\it tensor radical} $\cN$ of $\cT$ as
$$
\cN(X,Y):=\{f:X\pfeil Y\mid \|tr|fg=0\hbox{ for all }g:Y\pfeil X\}.
$$
One can show that $\cN$ is a tensor ideal (\cite{AKO}~7.1.1). If $K$
is a field then $\cN$ is the maximal proper tensor ideal of $\cT$
(\cite{AKO}~7.1.4.).

\Definition. An object $X$ of $\cT$ is called {\it
$\epsilon$\_semisimple} (or {\it $\epsilon$\_simple}) if
$\|End|_\cT(X)$ is a semisimple ring (or a division ring).

\Lemma Schur0. Let $S$ and $X$ be objects of $\cT$. Assume that $S$ is
$\epsilon$\_simple.
\Item{16} If $\cN(S,X)=0$ then every non\_zero morphism $S\pfeil X$
admits a retraction.
\Item{17} If $\cN(X,S)=0$ then every non\_zero morphism $X\pfeil S$
admits a section.\Par

\Proof: We prove \cite{I16}. The proof for \cite{I17} is
analogous. Let $f:S\pfeil X$ be a non\_zero morphism. Since
$f\not\in\cN(S,X)=0$ there is a morphism $g:X\pfeil S$ with
$\|tr|(gf)\ne0$. This implies that $gf$ is a non\_zero, hence
invertible endomorphism of $S$. Then $\tilde g:=(gf)^{-1}g$ is a
retraction of $f$.\qed

\noindent This implies the following Schur type lemma:

\Lemma Schur. Let $S_1$ and $S_2$ be two $\epsilon$\_simple objects
of $\cT$. Assume moreover that $\cN(S_1,S_2)=0$. Then every morphism
$S_1\pfeil S_2$ is either zero or an isomorphism.

\Proof: Let $f:S_1\pfeil S_2$ be non\_zero. By \cite{Schur0} there is
a morphism $g:S_2\pfeil S_1$, with $gf=1_{S_1}$. On the other hand,
$fg$ is a non\_zero idempotent, hence equal to $1_{S_2}$.\qed

\noindent $\epsilon$\_simple and $\epsilon$\_semisimple objects are
related in the following way:

\Proposition simple-semisimple. Let $\cT$ be a pseudo\_abelian tensor
category with $\cN=0$. Let $X$ be an object of $\cT$. Then the
following are equivalent:
\Item{20}$X$ is $\epsilon$\_semisimple.
\Item{21}$X$ is a direct sum of $\epsilon$\_simple objects.

\Proof: \cite{I20}$\Rightarrow$\cite{I21}: this direction works even
without the assumption $\cN=0$. By the structure theory of semisimple
rings we have
$$42
B:=\|End|_\cT(X)\cong M_{d_1}(K_1)\times\ldots\times M_{d_r}(K_r)
$$
where the $K_i$ are division rings. The canonical set of minimal
orthogonal idempotents of $B$ (their number is $\sum d_i$) splits $X$
as
$$41
X\cong X_1^{d_1}\oplus\ldots\oplus X_r^{d_r}
$$
with
$$
\|Hom|_\cT(X_i,X_j)=\cases{K_i&if $i=j$\cr0&if $i\ne j$\cr}
$$
The $X_i$ are, in particular, $\epsilon$\_simple.

\cite{I21}$\Rightarrow$\cite{I20}: Assume there is a decomposition
\cite{E41} such that $K_i:=\|End|_\cT(X_i)$ is a division ring and
such that $X_i\not\cong X_j$ for $i\ne j$. \cite{Schur} implies that
$\|Hom|_\cT(X_i,X_j)=0$ for $i\ne j$. This implies \cite{E42}.\qed

\Lemma dirsum. Let $\cT$ be a pseudo\_abelian tensor category with
$\cN=0$. Let $X_1$ and $X_2$ be two objects. Then $X_1\oplus X_2$ is
$\epsilon$\_semisimple if and only if both $X_1$ and $X_2$ are
$\epsilon$\_semisimple.

\Proof: If $X_1$ and $X_2$ are $\epsilon$\_semisimple then $X:=X_1\oplus
X_2$ is a direct sum of $\epsilon$\_simple objects. Thus, $X$ is
$\epsilon$\_semisimple. Assume conversely that $X$ is
$\epsilon$\_semisimple. Then the decomposition $X=X_1\oplus X_2$
corresponds to orthogonal idempotents $p_1$, $p_2$ of the semisimple ring
$B=\|End|_\cT(X)$. It is well\_known that $\|End|_\cT(X_i)=p_iBp_i$ is
again a semisimple ring.\qed

\noindent Here is our main criterion for semisimplicity:

\Corollary sscriterion. Let $\cT$ be a pseudo\_abelian tensor category
with $\cN=0$. Let $\cT'$ be a full subcategory which generates $\cT$
as a pseudo\_abelian category. Then $\cT$ is semisimple if and only if
every object of $\cT'$ is $\epsilon$\_semisimple.

The decompositions \cite{E42} and \cite{E41} are related in a more
canonical fashion which we recall now in a more general form. Let $B$
be a semisimple ring and let $\{M_\pi\mid\pi\in\hat B\}$ be a set
containing each simple $B$\_module up to isomorphism exactly
once. Then $K_\pi:=(\|End|_BM_\pi)^{\rm op}$ is a division ring and
$M_\pi$ is a $B\hbox{--}K_\pi$\_bimodule. Moreover,
$M_\pi^*:=\|Hom|_{K_\pi}(M_\pi,K_\pi)$ is a
$K_\pi\hbox{--}B$\_bimodule. With this notation, the decomposition
\cite{E42} corresponds to
$$23
B=\bigoplus_{\pi\in\Bh}\|End|_{K_\pi}M_\pi=
\bigoplus_{\pi\in\Bh}M_\pi\otimes_{K_\pi}M_\pi^*.
$$

Now assume an object $X$ of $\cT$ is endowed with a homomorphism
$B\pfeil\|End|_\cT(X)$. Then for any $\pi\in\Bh$ put
$$50
X^\pi:=M_\pi^*\otimes_BX.
$$
Here $V\otimes_BX$ is the object representing the functor
$Y\mapsto\|Hom|_B(V,\|Hom|_\cT(X,Y))$ (see, e.g., \cite{De}
formula~(3.7.1)). Then $X^\pi$ is a left $K_\pi$\_object of $\cT$. The
decomposition \cite{E41} becomes
$$13
X=\bigoplus_{\pi\in\Bh}M_\pi\otimes_{K_\pi}X^\pi
$$
Moreover, if $B=\|End|_\cT(X)$ then
$$44
\|Hom|_\cT(X^\pi,X^{\pi'})=\cases{K_\pi&if $\pi=\pi'$,\cr0&if $\pi\ne\pi'$.\cr}
$$

\beginsection core. The core of a relation

For general regular categories, it is difficult to control all
subobjects of a product $x\times y$. Therefore, in this and the next
section, we are going to restrict our attention to exact Mal'cev
categories (see \cite{Malcevdef}) because there all subobjects of a
product are basically pull\_backs.

More precisely, let $r\into x\times y$ be a relation. To get hold of
$r$ we first consider the images $\xq$ and $\yq$ of $r$ in $x$ and
$y$, respectively. Then we form the push\_out of $r\auf\xq$ along
$r\auf\yq$ (possible by \cite{pushpull}). Thus, we arrive at the
following diagram
$$3
\cxymatrix{
&&r\ar@{>>}[dl]\ar@{>>}[dr]&&\\
&\xq\ar@{ >.>}[dl]\ar@{.>>}[dr]&&\yq\ar@{.>>}[dl]\ar@{ >.>}[dr]\\
x&&c&&y
}
$$
where the square is a push\_out. The point is now, that in an exact
Mal'cev category $r$ can be recovered from the dotted part of the
diagram. In fact, \cite{pushpull} implies that the square is also a
pull\_back diagram. Thus, we obtain a bijection between subobjects of
$x\times y$ and isomorphisms between subquotients of $x$ and $y$ up to
some obvious equivalence. For the category of groups, this observation
is due to Goursat (\cite{Gou}~p.~47--48).

\Definition. Let $r\into x\times y$ be a relation. Then the object $c$
of diagram \cite{E3} is called the {\it core} of $r$.

\noindent The significance of this definition is summarized in the
following lemma.

\Lemma factorization. Let $\cA$ be an exact Mal'cev category, let
$r\into x\times y$ be a relation, and let $c$ be its core.
\Item{11} The morphism $\<r\>$ factorizes in $\cT(\cA,\delta)$ through $[c]$.
\Item{12} Assume that $\lambda\<r\>$, with $\lambda\ne0$ and
$r\ne\io$, factorizes in $\cT(\cA,\delta)$ through an object
$[z]$. Then $c\preceq z$.\Par

\Proof: \cite{I11} is obvious from diagram \cite{E3}.

\cite{I12} Assume $\lambda\<r\>$ is equal to the composition $[x]\pf G[z]\pf
H[y]$. Then $G$ and $H$ ``contain'' relations $s\into x\times z$ and
$t\into z\times y$ such that $r=t\circ s$. In other words, there is a
diagram
$$
\cxymatrix{
&&\tilde r\ar[ld]\ar[rd]&&\\
&s\ar[ld]\ar[rd]&&t\ar[ld]\ar[rd]&\\
x&&z&&y
}
$$
where the square is a pull\_back such that $r$ is the image of
$\tilde r$ in $x\times y$. Let $\zq$ be the image of $\tilde r$ in
$z$. Then we get the following diagram:
$$34
\cxymatrix{
&&\tilde r\ar@{=}[d]\ar@{>>}[dr]\ar@{>>}[dl]&&\\
&\sq\ar@{>>}[dl]\ar@{ >->}[d]\ar@{>>}[dr]&
\tilde r\ar[ld]|!{[l];[d]}\hole\ar[rd]|!{[r];[d]}\hole&
\tq\ar@{ >->}[d]\ar@{>>}[dl]\ar@{>>}[dr]&\\
\xq\ar @{ >->}[d]&
s\ar[ld]\ar[rd]&\zq\ar@{ >->}[d]&
t\ar[ld]\ar[rd]&\yq\ar@{ >->}[d]\\
x&&z&&y
}
$$
Here, the upper square is the pull\_back of the lower one by $\zq\into
z$ while $\xq$ and $\yq$ are the images of $\sq\pfeil x$ and
$\tq\pfeil y$, respectively. The morphisms $\sq\pfeil\zq$ and
$\tq\pfeil\zq$ are surjective since $\tilde r\pfeil\zq$ is. The upper
square is also a pull\_back diagram. Hence, the two morphisms from $\tilde
r$ to $\sq$ and $\tq$ are surjective, as well. This implies that $\xq$
and $\yq$ are the images of $\tilde r$ in $x$ and $y$, respectively.
By definition of $c$ there is a diagram
$$
\cxymatrix{
&\tilde r\ar@{>>}[ldd]\ar@{>>}[d]\ar@{>>}[rdd]\\
&r\ar@{>>}[ld]\ar@{>>}[rd]\\
\xq\ar@{>>}[rd]&&\yq\ar@{>>}[ld]\\
&c
}
$$
Since the upper square of \cite{E34} is also a push\_out
(\cite{pushpull}), we obtain a morphism $\zq\pfeil c$ which is
surjective since $\tilde r\pfeil c$ is. This yields the
desired subquotient diagram $\xymatrix@=10pt{z&\zq\ar@{
>->}[l]\ar@{>>}[r]&c}$.\qed

\noindent Here is the linearized version of the preceding theorem:

\Corollary decomposition. Let $x$ and $y$ be objects of an exact
Mal'cev category $\cA$.\Par
\Item{18} Every $\cT(\cA,\delta)$\_morphism $[x]\pfeil[y]$ factorizes
through an object of the form $[z_1]\oplus\ldots\oplus[z_n]$ with
$z_i\preceq x$ and $z_i\preceq y$ for all $i$.
\Item{19}Assume $x\not\prec x$ (see \cite{partialorder0}). Then there
is a decomposition
$$43
\|End|_{\cT}([x])=K[\|Aut|_{\cA}(x)]\oplus\|End|^\prec_{\cT}([x]).
$$
Here, $\|End|^\prec_{\cT}([x])$ is the two\_sided ideal of all
endomorphisms which factorize through an object of the form
$[z_1]\oplus\ldots\oplus[z_n]$ with $z_i\prec x$ for all $i$.\Par

\Proof: \cite{I18} Follows directly from
\cite{factorization}~\cite{I11}.

\cite{I19} Assume the core $c$ of a relation $r\into x\times x$ is not
a proper subquotient of $x$. Then the dotted arrows of diagram
\cite{E3} (with $y=x$) are all isomorphisms. Thus also the two solid
arrows are isomorphism which means that $r$ is the graph of an
automorphism of $x$. This and \cite{factorization}~\cite{I11} imply
that
$$
\|End|_{\cT}([x])=K[\|Aut|_{\cA}(x)]+\|End|^\prec_{\cT}([x]).
$$

To show that the sum is direct assume that the linear combination
$F=\sum_j\lambda_j[f_j]$ factorizes through
$[z_1]\oplus\ldots\oplus[z_n]$ with $z_i\prec x$ and with pairwise
different $f_j\in\|Aut|_\cA(x)$. Suppose $\lambda_j\ne0$. Then there
are relations $r\into x\times z_i$ and $s\into z_i\times x$ such that
the $\cT$\_composition $\<s\>\<r\>$ is a non\_zero multiple of
$[f_j]$. \cite{factorization}~\cite{I12} implies that
$x=\|core|(f_j)\preceq z_i$ in contradiction to $z_i\prec x$.\qed

\beginsection semisimplicity. The semisimplicity of $\cTq(\cA,\delta)$

We return to our pseudo\_abelian tensor category $\cT(\cA,\delta)$
attached to a regular category $\cA$ and a $K$\_valued degree function
$\delta$. In this section we address the problem whether
$$
\cTq(\cA,\delta):=\cT(\cA,\delta)/\cN 
$$
is a semisimple, hence abelian, tensor category. Except for very
degenerate cases (e.g., $\delta=1$), semisimplicity cannot be expected
unless all $\|Hom|$\_spaces are finite dimensional over
$K$. Therefore, we are going to assume that $\cA$ is {\it subobject finite},
i.e., that every object has only finitely many subobjects.

But even then there is a problem: Deligne (\cite{De}
Mise~en~garde~5.8) has constructed a pseudo\_abelian tensor category
over $\CC$ with finite dimensional $\|Hom|$\_spaces and $\cN=0$ which
is not semisimple.

To state our main criterion, let $\hat\cT$ be the class of of
isomorphism classes of pairs $(x,\pi)$ where $x$ is an object of $\cA$
and $\pi$ is an irreducible $K$\_representation of $\|Aut|_\cA(x)$.

\Theorem maintheorem. Let $\cA$ be a subobject finite, exact Mal'cev
category and let $\delta$ be a $K$\_valued degree function on $\cA$
where $K$ is a field of characteristic zero. Then:
\Item{22} $\cTq(\cA,\delta)$ is a semisimple (hence abelian) tensor
category.
\Item{23} For every $(x,\pi)\in\hat\cT$ there is, up to isomorphism,
at most one simple object $M$ of $\cTq(\cA,\delta)$ with
\itemitem{---}$M$ occurs in the isotypic component $[x]^\pi$ (notation
of \cite{E50}).
\itemitem{---}$M$ does not occur in $[y]$ for any $y\prec x$.
\Item{29}Let $\hat\cT_\delta$ be the set of $(x,\pi)\in\hat\cT$ such
that $M$ as in \cite{I23} exists and denote, in that case, $M$ by
$M^{(x,\pi)}$. Then $(x,\pi)\mapsto M^{(x,\pi)}$ is a bijection
between $\hat\cT_\delta$ and isomorphism classes of simple objects.
\Item{24} If $\cN=0$ then $\hat\cT_\delta=\hat\cT$.\Par

\Proof: In view of \cite{iddecomp} we may assume that $\*$ has no
proper subobject. This means, in particular, that
$K=\|End|_\cT(\1)$.

By \cite{sscriterion}, we have to show that every object of the form
$[x]$ is $\epsilon$\_semisimple, i.e., has a semisimple endomorphism
ring. By \cite{partialorder0} it suffices to prove the following
statement: {\it let $x$ be an object of $\cA$ such that $[y]$ is
$\epsilon$\_semisimple for all $y\prec x$. Then $x$ is
$\epsilon$\_semisimple.}

Let $\fS$ be the (finite) set of all $\epsilon$\_simple summands
occurring in some $[y]$ with $y\prec x$. If we apply \cite{Schur0}
successively to $[x]$ and all elements of $\fS$, we obtain a
decomposition
$$46
[x]=[x]_0\oplus[x]_1
$$
such that $[x]_1$ is a direct sum of elements of $\fS$ and 
$$
\|Hom|_\cT([x]_0,[y])=\|Hom|_\cT([y],[x]_0)=0
$$
for all $y\prec x$. The decomposition \cite{E34}
implies
$$
\|End|_{\cTq}([x])=K[\|Aut|_{\cA}(x)]+\|End|^\prec_{\cTq}([x])
$$
but the sum may no longer be direct. It is clear that $\|End|^\prec$ kills
$[x]_0$. Thus, we obtain a surjective homomorphism
$$45
K[\|Aut|_\cA(x)]\pfeil B:=\|End|_\cT([x]_0)
$$
Since $K$ is of characteristic zero and $\|Aut|_\cA(x)$ is a finite
group we conclude that $B$ is a semisimple ring. Thus $[x]_0$ and
therefore $[x]$ is $\epsilon$\_semisimple (\cite{dirsum}), showing
\cite{I22}.

Let
$$
[x]_0=\bigoplus_{\pi\in\hat B}M_\pi\otimes M^{(x,\pi)}
$$
be the $B$\_isotypic decomposition (see \cite{E13}). Then
$M^{(x,\pi)}$ is $\epsilon$\_simple (see \cite{E44}), hence
simple. Since $K[\|Aut|_\cA(x)]\pfeil B$ is surjective, we can think
of $(x,\pi)$ as being an element of $\hat\cT$.

If $M$ is any simple object as in \cite{I23} then $M$ can't appear in
$[x]_1$. Thus $M\cong M^{(x,\pi)}$ proving \cite{I23}.

The decomposition \cite{E46} shows, by induction, that $[x]$ is a
direct sum of objects $M^{(y,\pi)}$ with $y\preceq x$. In particular,
every simple object of $\cT$ is of the form $M^{(x,\pi)}$. Now assume
that there is an isomorphism $f:M^{(x,\pi)}\pf\sim
M^{(x',\pi')}$. This isomorphism extends to a morphism
$[x]\pfeil[x']$. \cite{decomposition}\cite{I18} implies that
$M^{(x,\pi)}$ occurs already in an object $[y]$ where $y$ is a
subquotient of both $x$ and $x'$. By definition of $M^{(x,\pi)}$ and
$M^{(x',\pi')}$, this subquotient can't be proper. Thus we obtain
$x\pf\sim x'$. We conclude $\pi=\pi'$ (see \cite{E44}), showing
\cite{I29}.

Finally \cite{I24} follows from the fact that \cite{E45} is an
isomorphism if $\cN=0$.\qed

\Example: [Deligne's case] Let $\cA=\FS^{\rm op}$, the opposite
category of the category of finite sets. Then $\hat\cT$ is the union
over $n\ge0$ of $\widehat S_n$. Therefore, $\hat\cT$ is parametrized
by Young diagrams of arbitrary size.

\Remark: I see no inherent reason why $\cTq(\cA,\delta)$ should only
be semisimple for exact Mal'cev categories. It is just the proof which
requires that condition. On the other hand, the description of simple
objects is probably only valid in the exact Mal'cev case.

\beginsection Moebius. The M\"obius algebra of a semilattice

In preparation for the next section, we review and refine in this
section some results from (semi-)lattice theory. Recall that a {\it
semilattice} is a set $L$ equipped with an associative and commutative
product $\wedge$ which is idempotent, i.e., with $u\wedge u=u$ for all
$u\in L$. Any semilattice is partially ordered by
$$
u\le v\Leftrightarrow u\wedge v=u.
$$
Conversely, one con recover the product from the partial order because
$u\wedge v$ is the largest lower bound of $\{u,v\}$.

If $L$ is any partially ordered set let $L^\io$ be $L$ with a new
minimum $\io$, i.e.,
$L^\io:=L\mathop{\cup\kern-4.6pt\raise5pt\hbox{.}\kern4pt}\{\io\}$
with $\io<u$ for all $u\in L$. We call $L$ a {\it partial semilattice}
if $L^\io$ is a semilattice. In analogy to (and, in fact, a special
case of) \cite{partreg}, a poset $L$ is a partial semilattice if and
only if any two element set $\{u,v\}\subseteq L$ either has no lower
bound at all or has an infimum.

\Examples: \Bsp. Let $x$ be an object of a category $\cA$. If $\cA$ is
complete regular then $\|sub|(x)$, the partially ordered set of
subobjects of $x$, is a semilattice. If $\cA$ is just regular then
$\|sub|(x)$ is a partial semilattice. That is why we are interested in
them.

\Bsp. If $L$ is a partial semilattice then every upper subset $U$ (i.e.,
   one with $u\in U, u\le v\Rightarrow v\in U$) is a partial
   semilattice.

\medskip

For a finite partial semilattice $L$ let $P$ be the $L\times
L$\_matrix with $P_{uv}=1$ if $u\le v$ and $P_{uv}=0$ otherwise. This
matrix is unitriangular and therefore has an inverse $M$. The entries
$\mu(u,v):=M_{uv}$ are the values of the {\it M\"obius function} $\mu$
of $L$. It is $\ZZ$\_valued with $\mu(u,v)=0$ unless $u\le v$.

The M\"obius function has a natural interpretation in terms of the
{\it M\"obius algebra} $A(L)$. If $L$ is a semilattice then $A(L)$ is
the free abelian group with basis $L$ and multiplication induced by
$\wedge$. For a partial semilattice, we put
$A(L)=A(L^\io)/\ZZ\io$. Thus $A(L)$ is the free abelian group over
$L$ and multiplication is induced by $\wedge$ with the proviso that
$u\wedge v=0$ if $u\wedge v$ does not exist in $L$.

Define the elements
$$
p_v=\sum_{u\le v}\mu(u,v)u\in A(L).
$$
Then one can show (see e.g. \cite{Greene}) that the $p_v$
form a basis of $A(L)$ with
$$32
p_u\wedge p_v=\delta_{u,v}p_v.
$$
and
$$
p_u\wedge v=\cases{p_u&if $u\le v$,\cr0&otherwise.\cr}
$$

Now we generalize a formula of Lindstr\"om, \cite{Li}, and
Wilf, \cite{Wi} to partial semilattices.

\Lemma Wilf. Let $\phi:L\pfeil K$ be a function on a finite
partial semilattice $L$. Then
$$33
\|det|\big(\phi(u\wedge v)\big)_{u,v\in L}=\prod_{w\in
L}\phi(p_w).
$$
where, on the left hand side $\phi(\io):=0$ while, on the right hand
side, $\phi$ is extended linearly to a map $A(L)\pfeil K$.

\Proof: The map $u\mapsto p_u$ is a unitriangular base change of
$A(L)$. Thus
$$
\|det|\big(\phi(u\wedge v)\big)_{u,v\in L}=
\|det|\big(\phi(p_u\wedge p_v)\big)_{u,v\in L}=
\prod_{w\in L}\phi(p_w).
$$\qed

Next, we need a generalization of a formula of Greene (\cite{Greene}
Thm.5) which compares the minimal idempotents $p_u$ of two M\"obius
algebras. For that we write $p^L_u$ and $\mu_L(u,v)$ to make the
dependence on $L$ explicit. Recall that a pair of maps $e^*:M\pfeil
L$, $e_*:L\pfeil M$ between posets is a {\it Galois connection} if
$$11
l\le e^*(m)\Longleftrightarrow e_*(l)\le m
$$
for all $l\in L$ and $m\in M$. In that case, it is known (see, e.g.,
\cite{EKMS} Prop.~3) that both maps are order preserving and that
$e^*$ preserves infima. In particular, if $L$ and $M$ are partial
semilattices then $e^*$ is multiplicative (with
$e^*(\io):=\io$).

\Lemma Stanley. Let $e^*:M\pfeil L$ and $e_*:L\pfeil M$ a Galois
connection between finite partial semilattices. For any
$l\in L$ put $m:=e_*(l)\in M$ and
$$
p^L_{l\pfeil m}:=\sum\limits_{l'\le l\atop
e_*(l')={m}}\mu_L(l',l)\ l'.
$$
Then
$$
p^L_l=e^*(p^M_m)\wedge p^L_{l\pfeil m}
$$

\Proof: Let $x:=e^*(p^M_m)\wedge p^L_l$. Then $x$ is a linear combination of
elements of the form $\tilde l:=e^*(m')\wedge l'$ with $m'\le m$ and
$l'\le l$. Clearly $\tilde l\le l$ and $\tilde l=l$ if $m'=m$ and
$l'=l$. Conversely, suppose $\tilde l=l$. Then $l'\ge l$ and therefore
$l'=l$. Moreover $e^*(m')\ge l$ implies $m'\ge e_*(l)=m$, hence
$m'=m$. This shows that $x=l+\hbox{lower order
terms}$. On the other hand, $x\in A(L)p^L_l=\ZZ p^L_l$ and
therefore $x=p^L_l$. Thus we have
$$
p^L_l=e^*(p^M_m)\wedge\sum_{l'\le l}\mu_L(l',l)\,l'
$$
We are done if we show that $l'\le l$ and $e^*(p^M_m)\wedge
l'\ne0$ implies $e_*(l')=m$. Put $m':=e_*(l')$. Then $e_*(l')\le m'$
implies $l'\le e^*(m')$, hence $l'=l'\wedge e^*(m')$. Thus also
$e^*(p^M_m)\wedge e^*(m')\ne0$. Hence $p^M_m\wedge m'\ne0$ which
implies $m\le m'$. On the other hand, $l'\le l$ implies $m'\le m$ and
therefore $m'=m$, as claimed.\qed

\beginsection determinant. The tensor radical of $\cT(\cA,\delta)$

As opposed to $\cTq(\cA,\delta)$, the categories $\cT(\cA,\delta)$ form
a nice family in dependence of $\delta$. In fact, assume for the
moment that $\cA$ is essentially small (i.e., equivalent to a small
category). Then we can define $K(\cA)$ as the commutative ring
generated by symbols $\<e\>$ (with $e$ surjective) and relations
\resetitem
\Item{}$\<1_x\>=1$ for all objects $x$,
\Item{}$\<\eq\>=\<e\>$ if $\eq$ is a pull\_back of $e$, and
\Item{}$\<e\eq\>=\<e\>\<\eq\>$ for all $e$, $\eq$ which can be
composed.

\noindent It is clear that
$$81
\Delta:\fE(\cA)\pfeil K(\cA):e\mapsto\<e\>
$$
is a universal degree function on $\cA$, i.e., every degree
function factorizes uniquely through $K(\cA)$. Moreover, the category
$\cT(\cA,\Delta)$ is a universal family of tensor categories in the
sense that
$$
\cT(\cA,\delta)={\rm idempotent\ closure\ of\
}\cT(\cA,\Delta)\otimes_{K(\cA)}K.
$$
Observe that all $\|Hom|$\_spaces of the universal category
$\cT(\cA,\Delta)$ are projective $K(\cA)$\_modules (of finite type,
in case $\cA$ is subobject finite).

Thus, since $\cT(\cA,\delta)$ is a ``fiber'' of the family
$\cT(\cA,\Delta)$ it is of interest when $\cT(\cA,\delta)$ itself is
semisimple, or at least when $\cT(\cA,\delta)=\cTq(\cA,\delta)$, i.e.,
when its tensor radical vanishes. The purpose of this section is to
give a simple numerical criterion.

\Definition. Let $K$ be a field. A $K$\_valued degree function
$\delta$ is {\it non\_singular} if the tensor radical $\cN$ of
$\cT(\cA,\delta)$ is $0$.

Assume in this section that $\cA$ is a subobject finite regular
category and that $K$ is a field (of any characteristic). Let $X$ be
an object of $\cT(\cA,\delta)$. If $\|sub|(\*)=\{\*\}$ then
$\|End|_\cT(\1)=K$ and there is a pairing
$$10
\beta_X:\|Hom|_{\cT}(\1,X)\times
\|Hom|_{\cT}(X,\1)\pfeil K:(G,F)\mapsto FG.
$$
Basically by definition, this pairing is non\_degenerate if and only
if $\cN(\1,X)=\cN(X,\1)=0$.

\Lemma pairing1. Assume $\|sub|(\*)=\{\*\}$. Then the pairing
$\beta_{[x]}$ is non\_degenerate if and only if
$$
\Omega_x:=\|det|\Big(\delta(u\cap v\pfeil\*)\Big)_{u,v\subseteq x}\ne0
$$
where we put $\delta(\io\pfeil\*)=0$.

\Proof: Every subobject $u$ of $x$ induces relations
$\delta_u:u\into\*\times x$ and $\epsilon_u:u\into x\times\*$. These
relations form a $K$\_basis in their respective
$\|Hom|_\cT$\_space. Moreover
$\beta_{[x]}(\delta_u,\epsilon_v)=\delta (u\cap v\pfeil\*)$. Thus
$\Omega_x$ is just the determinant of $\beta_{[x]}$ with respect
to these bases.\qed

We proceed with the calculation of $\Omega_x$. Observe that
$\|sub|(x)$ is a finite partial semilattice (see \S\cite{Moebius} for
the definition) with respect to intersection: $u\wedge v:=u\cap
v:=u\times_xv$. Thus, for any surjective morphism $e:x\auf y$ we may define
$$18
\omega_e:=\sum_{w\in\|sub|(x)\atop e(w)=y}\mu(w,x)\delta (w\auf y)\in K
$$
where $\mu$ is the M\"obius function of $\|sub|(x)$. If
$\|sub|(\*)=\{\*\}$ and $y=\*$ then this specializes to
$$
\omega_{x\auf\*}=\sum_{w\in\|sub|(x)}
\mu(w,x)\delta (w\auf\*).
$$
Now we have the following factorization:

\Lemma OMEGA. Assume $\|sub|(\*)=\{\*\}$. Then
$$
\Omega_x=\prod_{u\in\|sub|(x)}\omega_{u\auf\*}.
$$

\Proof: Apply \cite{Wilf} to $L=\|sub|(x)$ and
$\phi(u):=\delta(u\pfeil\*)$.\qed

The elements $\omega_{x\auf\*}$ factorize further. This is similar to
Stanley's factorization, \cite{St1}, of the characteristic polynomial
of a lattice.

\Lemma oMEGA. The element $\omega_e$ is multiplicative in $e$, i.e.,
if $x\Auf\eq y\Auf e z$ are surjective morphisms then
$\omega_{e\eq}=\omega_e\omega_\eq$.

\Proof: Every surjective morphism $e:x\auf y$
induces two maps
$$0
\eqalignno{&&e_*:\|sub|(x)\pfeil\|sub|(y):u\mapsto e(u)=
\|im|(u\into x\Auf ey)\cr
&&e^*:\|sub|(y)\pfeil\|sub|(x):v\mapsto x\times_yv\cr}
$$
They form a Galois connection since both $e_*(u)\le v$ and
$u\le e^*(v)$ are equivalent to $u\pfeil y$ factorizing through
$v$.

Now we apply \cite{Stanley} to $L=\|sub|(x)$, $M=\|sub|(y)$ and
$l=x$. Then $m=y$ since $e$ is surjective. Thus we get
$$12
\sum_{u\in\|sub|(x)}\mu(u,x)u=
\Bigg(\sum_{v\in\|sub|(y)}\mu(v,y)e^*(v)\Bigg)\cap
\Bigg(\sum_{w\in\|sub|(x)\atop e(w)=y}\mu(w,x)z\Bigg).
$$
Let $v\subseteq y$ and $w\subseteq x$ with $e(w)=y$. Put
$u:=e^*(v)\cap w$. Then the double pull\_back diagram
$$
\cxymatrix{
u\inj[r]\inj[d]&e^*(v)\inj[d]\sur[r]&v\inj[d]\\
w\inj[r]\ar@/_1pc/ @{>>}[rr]&x\sur[r]^e&y\sur[r]^\eq&z
}
$$
\smallskip\noindent
yields $e\eq(u)=z$ if and only if $\eq(v)=z$ and in this case
$$
\delta (u\auf z)=\delta (v\auf z)\cdot\delta (w\auf y).
$$
Thus, if we apply the function
$$
\phi(u):=\cases{\delta (u\pfeil z)&if $u\pfeil z$ is surjective\cr
0&otherwise\cr}
$$
to both sides of \cite{E12} then we get $\omega_{e\eq}=\omega_e\omega_\eq$.\qed

Let's call a surjective morphism $e$ {\it indecomposable} if in any
factorization $e=e'e''$ with $e',e''$ surjective precisely one of the
factors is an isomorphism. Combining the above yields:

\Corollary vanishing. Let $\cA$ be a subobject finite, regular
category with $\|sub|(\*)=\{\*\}$. Then the pairing
$\beta_{[x]}$ is non\_degenerate if and only if $\omega_e\ne0$ for all
indecomposable surjective morphisms $e:u\auf v$ with $u\preceq x$.

\Proof: The set $\|sub|(u)$ being finite for any object
$u$ implies that $u$ has only finitely many quotient objects. This, in
turn, entails that every surjective map is the composition of finitely
many indecomposable ones. We conclude with \cite{pairing1},
\cite{oMEGA}, and \cite{OMEGA}.\qed

From this we get our main vanishing theorem:

\Theorem. Let $\cA$ be a subobject finite, regular category and $K$ a
field. Then a $K$\_valued degree function $\delta$ is non\_singular
if and only if $\omega_e\ne0$ for all indecomposable surjective
morphisms $e$.

\Proof: \cite{iddecomp} reduces the assertion to the case
$\|sub|(\*)=\{\*\}$. If $\omega_{x\auf y}=0$ then $\Omega_x=0$ and
therefore $\cN(\1,[x])\ne0$. Conversely, the non\_vanishing of
$\omega_e$ for all indecomposable $e$ implies $\Omega_x\ne0$ and
therefore $\cN(\1,[x])=0$ for all $x$. From $\cN(\1,X\oplus
Y)=\cN(\1,X)\oplus\cN(\1,Y)$ we conclude $\cN(\1,X)=0$ for all objects
$X$ of $\cT$. This implies $\cN=0$ by \cite{E40}.\qed

\Corollary genericdegree. Let $\cA$ be an essentially small, subobject
finite, complete, exact, protomodular category. Then $\cA$ has
non\_singular degree functions.

\Proof: Consider the universal degree function $\Delta:\fE(\cA)\pfeil
K(\cA)$. \cite{degreecomplete}\cite{I75} implies that $K(\cA)$ is a
polynomial ring over $\ZZ$. For any surjective morphism $e:x\auf y$
consider $\omega_e$ computed with respect to $\Delta$. Then $\omega_e$
is a polynomial. \cite{protoDelta} asserts that the monomial
$\Delta(e)$ occurs only once in $\omega_e$ which implies
$\omega_e\ne0$. Now we can take for $K$ the field of fractions of
$K(\cA)$ (or any bigger field) and for $\delta$ the composition
$\fE(\cA)\pf\Delta K(\cA)\hookrightarrow K$.\qed

\Remark: We show in the last example below that there are exact
Mal'cev categories without non\_singular degree functions. So the
condition of protomodularity cannot be dropped.

\Examples: \Bsp. [Deligne's case] Let $\cA=\FS^{\rm op}$ where $\FS$
is the category of finite sets. Then surjective morphisms in $\cA$ are
injective maps in $\FS$. Let $t\in K$. The degree functions are
parametrized by $t\in K$ (see \cite{E57}). An injective map $e:A\into
B$ is indecomposable if $B\setminus e(A)=\{b\}$ is a one\_point
set. To compute $\omega_e$ we have to consider diagrams
$$
\cxymatrix{
A\inj[r]\inj[dr]&B\hbox to 0pt{$\,=A\cup\{b\}$\hss}\sur[d]\\
&Q
}
$$
There are two cases: either $Q=B$ or
$Q=A$. In the first case, $\mu(Q,B)=1$ and $\delta(A\into Q)=t$. The
second case depends on the image of $b$, so there are $|A|$
possibilities. Moreover, $\mu(Q,B)=-1$ and $\delta(A\into Q)=1$. This
implies $\omega_e=t-|A|$. Since $|A|$ is an arbitrary natural number
we conclude: {\it $\delta$ is non\_singular if and only if
$t\not\in\NN$.}

\Bsp. Let $G$ be a finite group and let $\cA$ be the opposite
category of the category of finite sets with a {\it free}
$G$\_action. This category is regular but not complete. Let
$\ell_0(A):=|A/G|$, the number of $G$\_orbits. Then all
degree functions are of the form
$$
\delta(e:A\into B)=t^{\ell_0(B\setminus e(A))}.
$$
Then as before one deduces that $e:A\into B$ is indecomposable if
$B\setminus e(A)$ is just one orbit and in that case
$\omega_e=t-|A|$. Thus {\it $\delta$ is non\_singular if and only if
$t\not\in\NN\,|G|$.}

\Bsp. Let $\cA$ be the category of non\_empty finite sets. Then there
is only one degree function with $\delta(e)=1$ for all surjective maps
$e$. The map $e:A\auf B$ is indecomposable if it identifies exactly
one pair $\{a_1,a_2\}$ of points to one point. There are exactly three
subsets $A_0\subseteq A$ such that $A_0\pfeil B$ is still
surjective. Hence $\omega_e=1-1-1=-1\ne0$. This shows that $\cN=0$ for
any field (even in positive characteristic). Because $\cA$ is not
Mal'cev we cannot apply \cite{maintheorem}. Thus, it is not clear
whether $\cT(\cA,\delta)$ is semisimple.

\Bsp. Let $\cA=\FV_{\FF_q}$ be the category of finite dimensional
$\FF_q$\_vector spaces. The degree functions are given by formula
\cite{E56}.The homomorphism $e$ is indecomposable if
$\|dim|\|ker|e=1$. To compute $\omega_e$ we have to consider diagrams
$$
\cxymatrix{
\hbox to 20pt{\hfill}&S\inj[d]\sur[dr]\\
&\hbox to 0pt{\hss$V\oplus\FF_q\cong\,$}U\sur[r]&V
}
$$
up to automorphisms of $S$. Again, there are two possibilities: $S=U$
or $S=V$. In the first case $\mu(S,U)=1$ and $\delta(S\auf V)=t$. In
the second case, $S$ is a section of $e$, hence there are $|V|$
possibilities. Since $\mu(S,U)=-1$ and $\delta(S\auf V)=1$ we get
$\omega_e=t-|V|$. We conclude: {\it $\delta$ is non\_singular if and
only if $t\not\in q^\NN$.} Observe that, in particular, $t=0$ is also
non\_singular.

\Bsp. Let, more generally, $\cA$ be a subobject finite abelian
category. The degree functions are given by formula \cite{E55}. A
surjective morphism $e:x\auf y$ is indecomposable if and only if
$s=\|ker|e$ is simple. A calculation as above shows that
$\omega_e=t_s-\alpha$ where $\alpha$ is the number of sections of
$e$. This number can be zero unless $s$ is injective. Otherwise,
$\alpha$ is a power of $q_s:=|\|End|_\cA(s)|$. We conclude: {\it
$\delta$ is non\_singular if and only if $t_s\not\in q_s^\NN$ for
all $s\in S$ and $t_s\ne0$ for all non\_injective $s\in S$.}

\Bsp. Let $\cA$ be the category of homomorphisms $f:U\pfeil V$ between
finite dimensional $\FF_q$\_vector spaces. This is the category of
$\FF_q$\_representation of the quiver $\bullet\pfeil\bullet$ and
therefore abelian. The simple objects are $s_1=(\FF_q\pfeil0)$ and
$s_2=(0\pfeil\FF_q)$. Only $s_1$ is injective. Let
$t_i:=t_{s_i}$. Then {\it $\delta$ is non\_singular if and only if
$t_1\not\in q^\NN$ and $t_2\not\in q^\NN\cup\{0\}$.}

\Bsp. Let $\cA$ be the category of (non\_empty) affine spaces over
$\FF_q$. The degree functions are given by
$$
\delta(e:X\auf Y)=t^{\|dim|X-\|dim|Y}.
$$
Moreover, $e$ is indecomposable if $X\cong Y\times \A^1$. In that
case, $\omega_e=t-\alpha$ where $\alpha$ is the number of section of
$e$. Since $\alpha=q\,|Y|$ we get {\it $\delta$ is non\_singular if
and only if $t\not\in q^{\NN^*}$}.

\Bsp. Let $\cA$ be the category of finite solvable groups. For a prime
$p$ let $v_p$ be the corresponding valuation of $\ZZ$. Then all
degree functions on $\cA$ are given by
$$
\|deg|(e:G\auf H)=\prod_pt_p^{v_p(|\|ker|e|)}
$$
with infinitely many parameters $t_2,t_3,t_5,\ldots\in K$. The map
$e:G\auf H$ is indecomposable if and only if $K=\|ker|e$ is a minimal
non\_trivial normal subgroup of $G$. Then $K$ is an elementary abelian
group of order $p^n$, say. Let $L\subseteq G$ be a subgroup with
$e(L)=H$, i.e., $G=KL$. The intersection $L\cap K$ is a subgroup of
$G$ which is normalized by $K$ (since $K$ is abelian) and $L$ hence by
$G=KL$. Minimality of $K$ implies either $K\subseteq L$ or $K\cap
L=1$. In the first case, we have $L=G$, $\mu(L,G)=1$ and $\|deg|(L\auf
H)=t_p^n$. In the second case holds $G=L\semidir K$. Thus
$\mu(L,G)=-1$ and $\|deg|(L\auf H)=1$. The number of complements can
be zero in which case $\omega_e=t_p^n$. Otherwise, the conjugacy
classes of complements are parametrized by $H^1(H,K)$ and each
conjugacy class has $|K/K^H|$ elements. Thus, the number of
complements is a power $p^N$ with $N\in\NN$. We conclude that {\it
$\delta$ is non\_singular if and only if $t_p\ne0$ and $t_p\ne p^r$
for all primes $p$ and all $r\in\QQ_{\ge0}$.}

\Bsp. We present an example for which {\it every} degree function is
degenerate. Let $\cA$ be the category of finite pointed Mal'cev
algebras. Objects of this category are finite sets $A$ equipped with a
base point $0\in A$ and a ternary operation $m:A^3\pfeil A$ satisfying
the identities $m(a,a,c)=c$ and $m(a,c,c)=a$ for all $a,c\in A$. Now
we define ternary operations $m_3$ and $m_2$ on $A:=\{0,1,2\}$ and
$B:=\{0,1\}$, respectively by
$$58
m_3(a_1,a_2,a_3)=\cases{
a_3&if $a_1=a_2$\cr
a_1&if $a_2=a_3$\cr
0&if exactly one of $a_1,a_2,a_3$ equals $0$\cr
1&otherwise\cr
}
$$
and
$$
m_2(b_1,b_2,b_3)=b_1-b_2+b_3\ \|mod| 2
$$
Then one verifies easily:
\item{$\bullet$}$m_3$ and $m_2$ are Mal'cev operations.
\item{$\bullet$}The only (pointed) subalgebras of $A$ are $A_0=\{0\}$,
$A_1=\{0,1\}$, and $A$ (the crucial point is $m_3(0,2,0)=1$). The only
subalgebras of $B$ are $B_0=\{0\}$ and $B$.
\item{$\bullet$}The map $e:A\auf B$ with $e(0)=0$ and $e(1)=e(2)=1$ is
an $\cA$\_morphism. For this notice that in the ``otherwise'' case of
\cite{E58} either none or exactly two of $a_1,a_2,a_3$ are equal to $0$.\Par

\noindent Since $A_0=\phi^{-1}(B_0)$, the pull\_back of $e$ by $B_0\pfeil B$ is
an isomorphism. Hence $\delta(e)=1$ and
$$
\omega_e=\mu(A,A)\delta(A\auf B)+\mu(A_1,A)\delta(A_1\auf
B)=1\cdot1+(-1)\cdot1=0.
$$
Thus, {\it $\delta$ is always degenerate}. Observe that $\cA$ is a
pointed, exact Mal'cev category.

\medskip

If one traces through the proof of \cite{maintheorem} and analyzes
exactly which $\cN(X,Y)$ have to be zero, one obtains the
following semisimplicity statement:

\Theorem. Let $\cA$ be a subobject finite, exact Mal'cev category, $K$
a field of characteristic zero and $\delta$ a $K$\_valued degree
function. Assume $\omega_f\ne 0$ for every indecomposable $f:u\auf v$
with $u\preceq x\times y$ for some $y\prec x$. Then $\|End|_\cT([x])$
is a semisimple $K$\_algebra.

\Example: If $\cA=\FS^{\rm op}$ then $\|End|_\cT([x])$ is known as {\it
partition algebra} (Martin \cite{Martin}). If $x=A$ is a set with $n$
elements then $y$ and $x\times y$ have at most $n-1$ and $2n-1$
elements, respectively. Thus $|v|\le2n-2$. We conclude that
$\|End|_\cT([x])$ is semisimple for $t\ne0,\ldots,2n-2$ which was
first proved by Martin\_Saleur \cite{MS}. Similarly, the
$\FF_q$\_analog of the partition algebra is semisimple unless
$t=1,q,\ldots,q^{2n-2}$.

\beginsection Tannakian. Tannakian degree functions

In this section, we investigate tensor functors from $\cT(\cA,\delta)$
to the category of $K$\_vector spaces. This will answer in particular
when $\cTq(\cA,\delta)$ is Tannakian, i.e., equivalent to
$\|Rep|(G,K)$ for some pro\_algebraic group over $K$, at least if $K$
is algebraically closed of characteristic zero.

Let $\FS$ be the category of {\it finite} sets and let $\FV_K$ be the
tensor category of finite\_dimensional $K$\_vector spaces. There is a
functor $\FS\pfeil\FV_K$ which maps a set $A$ to $K[A]$, the vector
space with basis $A$. Let $\pi:A\pfeil B$ be a map. Then
$K[\pi]$, also denoted by $\pi$, is the homomorphism
$$
\pi:K[A]\pfeil K[B]:a\mapsto\pi(a).
$$
Its transpose is the homomorphism
$$
\pi^\vee:K[B]\pfeil K[A]:b\mapsto \sum_{a\in\pi^{-1}(b)}a.
$$

\Lemma squarediagram. Consider the commutative diagram of sets
$$51
\cxymatrix{
A\ar[r]^{\pi_1}\ar[d]_{\tau_1}&B\ar[d]^{\tau_2}\\
C\ar[r]_{\pi_2}&D
}
$$
If \cite{E51} is a pull\_back then
$$20
\tau_2^\vee\pi_2=\pi_1\tau_1^\vee:\ K[C]\pfeil K[B].
$$
The converse is true if $\|char|K=0$.

\Proof: For $c\in C$ let $A_c:=\tau_1^{-1}(c)$, $d:=\pi_2(c)$, and
$B_d:=\tau_2^{-1}(d)$. Diagram \cite{E51} is a pull\_back if and
only if the map $\pi_c:A_c\pfeil B_d$ induced by $\pi_1$ is an isomorphism
for all $c\in C$. Now the assertion follows from
$$
\eqalign{
\tau_2^\vee\pi_2(c)&\,=\sum_{b\in B_d}b,\cr
\pi_1\tau_1^\vee(c)&\,=\sum_{b\in B_d}|\pi_c^{-1}(b)|b.\cr
}
$$\qed

A map $\pi$ between two finite sets is called {\it uniform} if all of
its fibers have the same cardinality. If $\pi:A\pfeil B$ is uniform
and $B\ne\emptyset$ then we call $\|deg|\pi:=|A|/|B|$ the {\it degree
of $\pi$}. In other words, $\|deg|\pi$ is the cardinality of the
fibers of $\pi$. Observe that $\emptyset\pfeil\emptyset$ has no
degree.

Let $\cA$ be a regular category.

\Definition unidef. Let $P:\cA\pfeil\FS$ be a functor and extend it to a
functor $P^\io:\cA^\io\pfeil\FS$ by setting $P^\io(\io)=\emptyset$.
Then $P$ is called {\it uniform} if
\item{---}$P^\io$ preserves finite limits (i.e., is {\it left exact}),
and
\item{---}$P$ maps surjective morphisms to uniform maps.

\Remark: In terms of $P$ alone, the first condition says that $P$
preserves finite limits and that $P(u)\times_{P(y)}P(v)=\emptyset$
whenever $u\times_yv$ does not exist. In particular, for a complete
regular category the first condition could be replaced by ``$P$ left
exact''.

\Definition adapteddef. A degree function $\delta:\fE(\cA)\pfeil K$ is {\it
adapted} to a uniform functor $P:\cA\pfeil\FS$ if
$$48
\delta(x\auf y)=\|deg|(P(x)\pfeil P(y))\ {\rm whenever}\
P(y)\ne\emptyset.
$$

\Remark: Call a uniform functor $P$ {\it non\_degenerate} if
$P(x)\ne\emptyset$ for all $x$. In that case, it is easy to check that
\cite{E48} defines a degree function on $\cA$. Thus, for a
non\_degenerate uniform functor $P$ there is precisely one degree
function $\delta_P$ adapted to it. Observe that if $\cA$ is pointed
then all uniform functors are non\_degenerate. Indeed, the
left\_exactness of $P$ implies that $P(\*)$ is a terminal object of
$\FS$, i.e., a one\_point set. Since $\cA$ is pointed, the unique map
$\*\pfeil x$ induces $P(\*)\pfeil P(x)$ which implies
$P(x)\ne\emptyset$.

\Example: [Deligne's case] At this point, we give only one
example. There will be more after \cite{UniThm2}. If $\cA=\FS^{\rm
op}$ and $X$ is any finite set then $P(A):=\|Hom|_\FS(A,X)$ is uniform
functor. In fact, it is clearly left exact. Moreover, let $e:A\into B$
be injective then any map $f:A\pfeil X$ can be extended to $B$ by
freely choosing the extension on $B\setminus e(A)$. Thus, the number
of extensions is $|X|^{|B\setminus e(A)|}$, independent of $f$. There
is exactly one degree function adapted to $P$, namely the one with
$t=|X|$ (notation of \cite{E57}).

\medskip

\Theorem. Let $\cA$ be a regular category and $K$ a
field.\Par
\Item{25}Let $P:\cA\pfeil\FS$ be a uniform functor and $\delta$ a
$K$\_valued degree function. Then there is a
tensor functor $T_P$ such that the diagram
$$22
\cxymatrix{
\cA\ar[r]^P\ar[d]_{[*]}&\FS\ar[d]^{K[*]}\\
\cT(\cA,\delta)\ar[r]^(.55){T_P}&\FV_K
}
$$
commutes if and only if $\delta$ is adapted to $P$. Moreover, $T_P$ is
unique.\Par
\Item{26}Assume that $K$ is algebraically closed of characteristic
zero. Let $T:\cT(\cA,\delta)\pfeil\FV_K$ be a tensor functor. Then
there is uniform functor $P:\cA\pfeil\FS$ (unique up to equivalence)
such that $\delta$ is adapted to $P$ and $T$ is equivalent to
$T_P$.\Par

\Proof: \cite{I25} Assume $\delta$ is adapted to $P$. Because
$\cT(\cA,\delta)$ is the universal pseudo\_abelian extension of
$\cT^0(\cA,\delta)$, it suffices to construct $T_P$ on $\cT^0$. The
commutativity of \cite{E22} forces us to define on objects
$$25
T_P([x])=K[P(x)].
$$

If $r\pfeil x\times y$ is a relation then, as a $\cT$\_morphism,
$\<r\>=[r\pfeil y][r\pfeil x]^\vee$. Thus, we have to define on morphisms
$$26
T_P(\<r\>)=P(r\pfeil y)P(r\pfeil x)^\vee:K[P(x)]\pfeil K[P(y)].
$$
This shows uniqueness.

Next, we show that \cite{E25} and \cite{E26} define indeed a tensor
functor. For that, consider the diagrams
$$28
\cxymatrix{
&&t\ar[ld]\ar[rd]\\
&r\ar[ld]\ar[rd]&&s\ar[ld]\ar[rd]\\
x&&y&&z
}
\qquad
\cxymatrix{
&t\sur[d]^e\\
&\tq\ar[ld]\ar[rd]\\
x&&z
}
$$
where the square is a pull\_back (in $\cA^\io$) and where $\tq$ is the
image of $t$ in $x\times z$. Apply $P$ to \cite{E28}. Then the
left\_exactness of $P^\io$ and \cite{E20} imply
$$
\eqalign{
T_P(\<s\>)T_P(\<r\>)&=
P(s\pfeil z)P(s\pfeil y)^\vee P(r\pfeil y)P(r\pfeil x)^\vee=\cr
&=P(s\pfeil z)P(t\pfeil s)P(t\pfeil r)^\vee P(r\pfeil x)^\vee=\cr
&=P(t\pfeil z)P(t\pfeil x)^\vee=
P(\tq\pfeil z)P(e)P(e)^\vee P(\tq\pfeil x)^\vee.\cr
}
$$
This is zero if $P(\tq)=\emptyset$. Otherwise
the uniformity of $P(e)$ implies
$$
P(e)P(e)^\vee=\|deg|P(e)\cdot1_{P(\tq)}
$$
and therefore
$$
T_P(\<s\>)T_P(\<r\>)=\|deg|P(e)\cdot T_P(\<\tq\>).
$$
On the other hand,
$$
T_P(\<s\>\<r\>)=T_P(\delta(e)\<\tq\>)=\delta(e)\cdot T_P(\<\tq\>).
$$
Thus, $T_P$ is a functor if and only if $\delta$ is adapted
to $P$. The fact that $T_P$ is a tensor functor follows from
$$
\eqalign{T_P([x]\otimes[y])&=T_P([x\times y])=K[P^\io(x\times y)]=\cr
&=K[P(x)\times P(y)]=K[P(x)]\otimes K[P(y)]=T_P(x)\otimes T_P(y)\cr}
$$

Now we prove part \cite{I26} of the theorem. The map $x\mapsto
(x):=[x]^\vee$ defines a {\it contravariant} functor from $\cA$ to
$\cT(\cA,\delta)$. It still has the property $(x\times y)=(x)\otimes
(y)$. Thus, the unique morphism $x\pfeil\*$ and the
diagonal morphism $x\into x\times x$ define $\cT$\_morphisms
$$
(x)\pfeil\1\quad{\rm and}\quad (x)\otimes(x)\pfeil (x)
$$
which equip $(x)$ with the structure of a unital commutative ring
object. Every multiplication $m:X\otimes X\pfeil X$ induces a trace
on $X$,
$$
\|tr|:X=X\otimes\1\Pf{1_X\otimes\|ev|^\vee}X\otimes X\otimes X^\vee
\Pf{m\otimes1_{X^\vee}}X\otimes X^\vee\Pf\sim X^\vee\otimes X\Pf{\|ev|}\1,
$$
and therefore a trace form
$$
X\otimes X\Pf mX\Pf{\|tr|}\1.
$$
An easy calculation shows that the trace on $(x)$ is
induced by the relation $x\pf\sim\*\times x$. Therefore, the trace
form comes from the relation
$$
\cxymatrix{
&x\ar[dl]\ar[dr]^{\rm diag}\\
\*&&x\times x
}
$$
This relation is just the evaluation morphism of the selfduality
$(x)^\vee=(x)$ (see \cite{E47}).

Let now $T:\cT(\cA,\delta)\pfeil\FV_K$ be a tensor functor. Then
$\cO(x):=T\big((x)\big)=T([x]^\vee)$ inherits all the properties above from
$(x)$: it is a finite\_dimensional unital commutative $K$\_algebra
such that the trace form is non\_degenerate. Since $K$ is
algebraically closed, this implies that $\cO(x)$ is isomorphic to
$K\times\ldots\times K$ as a $K$\_algebra. More canonically, put
$$
P(x):=\|Spec|\cO(x)=\|AlgHom|_K(\cO(x),K).
$$
Then $P$ is a covariant functor $\cA\pfeil\FS$ such that the functors
$T_P(x)=K[P(x)]$ and $T([x])$ from $\cA$ to $\FV_K$ are equivalent.

The proof of uniqueness above did not use any properties of $P$. This
implies that also $T_P$ and $T$ are equivalent to each other. It
remains to show that $P$ is uniform and that $\delta$ is adapted to
$P$. For that, consider the following sequence of diagrams in $\cA$
$$
\cxymatrix{
&&t\ar[dl]_\gq\ar[dr]^\fQ\\
&r\ar@{.>}[dl]_\sim\ar[dr]_f&&s\ar[dl]^g\ar@{.>}[dr]^\sim\\
r&&y&&s
}
\longleftrightarrow
\cxymatrix{
&t\ar[dl]_\gq\ar[dr]^\fQ\\
r&&s
}
\longleftrightarrow
\cxymatrix{
&&t\ar[dl]_\sim\ar[dr]^\sim\\
&t\ar[dl]_\gq\ar[dr]_\sim&&t\ar[dl]^\sim\ar[dr]^\fQ\\
r&&t&&s
}
$$
where the left square is an arbitrary pull\_back. It implies the
relation $[g]^\vee[f]=[\fQ][\gq]^\vee$. Thus we also have $P(g)^\vee
P(f)=P(\fQ)P(\gq)^\vee$. From \cite{squarediagram} we infer that also
$$
\cxymatrix{
&&P(t)\ar[dl]_{P(\gq)}\ar[dr]^{P(\fQ)}\\
&P(r)\ar[dr]_{P(f)}&&P(s)\ar[dl]^{P(g)}\\
&&P(y)&&
}
$$
is a pull\_back. Since $P(\io)$ and $P(\*)$ have to be the empty and
one\_point set, respectively, we see that $P$ is left exact.

Finally, let $e:x\auf y$ be surjective, inducing an algebra morphism
$(e):(y)\pfeil(x)$. The restriction of $\|tr|_{(x)}:(x)\pfeil\1$ to
$(y)$ is given by the relation $x\pfeil\*\times y$ which implies
$\|tr|_{(x)}\circ(e)=\delta(e)\|tr|_{(y)}$. Thus also for
$\cO(y)\pfeil\cO(x)$ holds that the restriction of $\|tr|_{\cO(x)}$ to
$\cO(y)$ is $\delta(e)\|tr|_{\cO(y)}$. But that means that the map on
spectra $P(x)\pfeil P(y)$ is uniform of degree $\delta(e)$ (if
$P(y)\ne\emptyset$). This shows that $\delta$ is adapted to $P$.\qed

Finally, we investigate kernel and image of a tensor functor
$T:\cT(\cA,\delta)\pfeil\FV_K$. Using Tannaka theory (see
e.g. \cite{Saa}) this would be easy if $\cT(\cA,\delta)$ were an
abelian category or, more generally, if $T$ would factorize through
$\cTq(\cA,\delta)$. This is actually true but will be shown only a
posteriori.

For any profinite group $G$ let $\FS(G)$ be the category of finite
sets equipped with a continuous $G$\_action. Let $\|Rep|(G,K)$ be the
category of continuous {\it finite dimensional} representations of $G$
over $K$. Then $A\mapsto K[A]$ provides also a functor
$\FS(G)\pfeil\|Rep|(G,K)$.

Let $P:\cA\pfeil\FS$ be a uniform functor. If $\cA$ is essentially
small then $\|Aut|P$ is an example of a profinite group. Indeed,
$\|Aut|P$ is the closed subgroup of $\prod_xS_{P(x)}$ defined by the
equations $P(f)\pi_x=\pi_yP(f)$. Here, $S_{P(x)}$ is the symmetric
group on a set $P(x)$. Moreover, $x$ and $f:x\pfeil y$ run through all
objects and morphisms of a skeleton of $\cA$. Clearly $P$ will factorize
canonically through $\FS(\|Aut|P)$ and diagram \cite{E22} can be
refined to
$$
\cxymatrix{
\cA\ar[r]^P\ar[d]&\FS(\|Aut|P)\ar[d]\\
\cT(\cA,\delta)\ar[r]^(.4){\cT_P}&\|Rep|(\|Aut|P,K)
}
$$
(provided $\delta$ is adapted to $P$).

For every injective morphism $y\into x$ the map
$P(y)\pfeil P(x)$ is injective, as well. Thus we may define
$$
P^*(x):=P(x)\setminus\bigcup_{y\subsetneq x}P(y).
$$
Note that $P^*$ is functorial for surjective morphisms only. Since $P^\io$
preserves intersections, for every $a\in P(x)$ there is a unique
minimal subobject $z$ with $a\in P(y)$. This means that $P(x)$ is the
disjoint union of the sets $P^*(y)$ with $y\in\|sub|(x)$.

From now on we assume that $\cA$ is an essentially small, subobject
finite, regular category and that $P:\cA\pfeil\FS$ is a uniform functor.

\Lemma. Let $e:x\auf y$ be a surjective morphism. Then $P^*(x)\pfeil
P^*(y)$ is uniform.

\Proof: For $a\in P^*(y)$ let $P(x)_a\subseteq P(x)$ and
$P^*(x)_a\subseteq P^*(x)$ be the fibers over $a$. Since $P(x)_a$ is
the disjoint union of $P^*(z)_a$ with $z\subseteq x$, M\"obius inversion
yields
$$
|P^*(x)_a|=\sum_{z\subseteq x}\mu(z,x)|P(z)_a|.
$$
Because $a\in P^*(y)$, it suffices to sum over those $z$ with $z\pfeil
y$ surjective. In that case, uniformity of $P$ implies that $|P(z)_a|$ is
independent of $a$. Hence the same holds for $|P^*(x)_a|$.\qed

The next assertion shows that $\|Aut|P$ is sufficiently big.

\Lemma uniform2. Let $x$ be an object of $\cA$. Then $\|Aut|P$ acts
transitively on $P^*(x)$.

\Proof: According to \cite{AM} Cor.~App.~2.8, every left exact functor
from $\cA$ to the category of (not necessarily finite) sets is
pro\_representable. This means that there is a small filtering
category $\cC$ and a functor $\cC^{\rm op}\pfeil\cA: i\mapsto p_i$
such that $P(u)=\|lim|\limits_{\longrightarrow}\|Hom|_\cA(p_i,u)$ for
all objects $u$ of $\cA$. We denote the corresponding pro-object by
$p$. Then $P^*(x)$ is the set of surjective morphisms $p\auf x$.

So, let $a,b:p\auf x$ be two surjective morphisms. We have to show
that there is an automorphism $f$ of $p$ such that $a=bf$.
Our plan is to construct $f$ by lifting $b:p\pfeil x$ to a morphism
$p\pfeil p$ as in the following diagram:
$$
\cxymatrix{
p\sur[ddrr]_b\ar@{-->}[rr]^f\ar@{.>>}[rrd]^{f_j}&&p\sur[d]\sur@/^2pc/[dd]^a\\
&&p_j\sur[d]^{a_j}\\
&&x
}
$$
For that observe that $a:p\auf x$ is represented by a surjective
morphism $p_i\pfeil x$ for some object $i$ of $\cC$. By replacing
$\cC$ with a cofinal subcategory, we may assume that $i$ is an initial
object of $\cC$. Moreover, the descending chain condition for
subobjects implies that we may assume that all morphisms $p_j\pfeil
p_k$ are surjective (simply replace $p_j$ by the intersection of all
images of the $p_i\pfeil p_j$).

Now a surjective morphism $f_j:p\auf p_j$ corresponds to an element of
$P^*(p_j)$. By \cite{uniform2}, the $P^*(p_j)$ form a projective system
of finite sets such that all structure maps $P^*(p_k)\pfeil P^*(p_j)$
are uniform. The canonical projection $p\pfeil p_j$ furnishes an element of
$P^*(p_j)$ which shows that this set is not empty. We conclude that
all maps $P^*(p_k)\pfeil P^*(p_j)$ are surjective. Now let
$P^*(p_j)_b$ be the fiber over $b\in P^*(x)$. Then also the
$P^*(p_j)_b$ form a projective system of finite sets such that all
structure maps are surjective. It follows that
$\|lim|\limits_{\longleftarrow}P^*(p_j)_b$ is non\_empty. Any element
of that limit corresponds to an endomorphism $f$ of $p$ such that $b=af$.

This $f$ induces an endomorphism $\phi$ of the functor $P$ with
$\phi(a)=b$. It remains to show that $\phi$ is invertible. By
construction, $f$ is represented by surjective morphisms $p\auf
p_j$. This implies that for any object $y$ of $\cA$ the self\_map
$\phi_y:P(y)\pfeil P(y)$ is injective. But then it is invertible since
$P(y)$ is finite.\qed

\Lemma. Assume $\delta$ is adapted to $P$. Then the functor
$$
T_P:\cT(\cA,\delta)\pfeil\|Rep|(\|Aut|P,K)
$$
is full, i.e., surjective on $\|Hom|$\_spaces.

\Proof: Let $X$, $Y$ be objects of $\cT(\cA,\delta)$. Because of
$\|Hom|_\cT(X,Y)=\|Hom|_\cT(\1,X^\vee\otimes Y)$ and because $T_P$
is a tensor functor we may assume $X=\1$. Moreover, since $T_P$
commutes with direct sums, is suffices to consider $Y=[x]$. Thus we
have to show that
$$39
\|Hom|_\cT(\1,[x])\Pfeil P(x)^G
$$
is surjective for all objects $x$ of $\cA$ (with $G:=\|Aut|P$). By
\cite{uniform2}, the $G$\_orbits in $P(x)$ are precisely the
non\_empty sets among the sets $P^*(y)$ with $y\in\|sub|(x)$. For a
subset $A$ of $P(x)$ let $\sum U=\sum_{\alpha\in A}\alpha$. Thus the
right hand space of \cite{E39} is spanned by all $\sum P^*(y)$ with
$y\in\|sub|(x)$ (it is not a basis since some of these sums may be
zero). By triangularity, $P(x)^G$ is also spanned by the elements
$\sum P(y)$ with $y\in\|sub|(x)$. But $\sum P(y)=T_P(F)$ where $F$
is the relation
$$
\cxymatrix{
&y\ar[dl]\inj[dr]\\
\*&&x.
}
$$\qed

Now we can prove the specialization theorem:

\Theorem UniThm. Let $\cA$ be an essentially small, subobject finite,
regular category, let $P:\cA\pfeil\FS$ be a uniform functor, and let
$\delta$ be $K$\_valued degree function adapted to $P$ where $K$ is a
field of characteristic zero. Then $T_P$ induces an
equivalence of tensor categories
$$29
\cTq(\cA,\delta)\Pf\sim\|Rep|(\|Aut|P,K).
$$

\Proof: Since $T_P$ is full and preserves traces, a $\cT$\_morphism
is in $\cN$ if and and only if its image is. Because $K$ is of
characteristic zero, the category $\|Rep|(\|Aut|P,K)$ is
semisimple. This shows that $T_P$ factorizes through
$\cTq(\cA,\delta)$ and that the functor \cite{E29} fully faithful. The action
of $G$ on the entirely of all sets $P([x])$ is effective. Therefore,
$\|Rep|(\|Aut|P,K)$ is generated, as a pseudo\_abelian tensor
category, by the representations of the form $T_P([x])$. This
implies that \cite{E29} is an equivalence.\qed

\Corollary UniThm2. Let $K$ be an algebraically closed field of
characteristic zero. Then $\cTq(\cA,\delta)$ is Tannakian if and only
if $\delta$ is adapted to some uniform functor $P:\cA\pfeil\FS$.

Before we go on with examples we introduce the following language:

\Definition. Let $\cA$ be subobject finite, regular category, $A$ a
commutative domain with field of fractions $F$, and $\Delta$ a
$A$\_valued degree function on $\cA$. Let $\fT:=\{\cT_i\mid i\in I\}$
be a family of semisimple tensor categories such that
$K_i:=\|End|_{\cT_i}(\1)$ is a field for all $i\in I$. Then we say
that {\it $\cT(\cA,\Delta)$ interpolates the family $\fT$} if:
\Item{76}The $F$\_valued degree function of $\cA$ induced by
$\Delta$ is non\_singular.
\Item{}There is a family of homomorphism $\phi_i:A\pfeil K_i$ such that
\itemitem{a)}the product homomorphism
$\prod_i\phi_i:A\pfeil\prod_iK_i$ is injective and
\itemitem{b)}the category $\cT_i$ is tensor equivalent to
$\cTq(\cA,\phi_i\circ\Delta)$ for all $i\in I$.

\noindent Observe that \cite{I76} is equivalent to $\omega_e\ne0$ for
all indecomposable surjective morphisms $e$ where we consider
$\omega_e$ as an element of $A$. By \cite{genericdegree} this is
automatically true for protomodular categories.

Every object of $\cT(\cA,\Delta)$ can be
considered as an object of $\cTq(\cA,\phi_i\circ\Delta)$ and, via the
equivalence in b), as an object of $\cT_i$. For any two objects $X$ and
$Y$ we get maps
$$
\|Hom|_{\cT(\cA,\Delta)}(X,Y)\otimes_AK_i\Auf\Phi
\|Hom|_{\cTq(\cA,\phi_i\circ\Delta)}(X,Y)\Pf\sim
\|Hom|_{\cT_i}(X,Y)
$$
The homomorphism $\Phi$ is always surjective and its kernel is it the
tensor radical of $\cT(\cA,\phi\circ\Delta)$. Its vanishing is
equivalent to the non\_vanishing (in $K_i$) of certain $\omega_e$
which are finite in number. Condition a) says that the $K_i$\_valued
points $\phi_i$ of $\|Spec|A$ are Zariski dense. Together this shows
that also $\Phi$ is an isomorphism for infinitely many $i\in I$. If,
for example, $\|Hom|_{\cT(\cA,\Delta)}(X,Y)$ were a free $A$\_module
then we would get a basis of all infinitely many of the spaces
$\|Hom|_{\cT_i}(X,Y)$. Moreover, the structure constants of the
composition
$$
\|Hom|_{\cT_i}(X,Y)\times \|Hom|_{\cT_i}(Y,Z)\Pfeil\|Hom|_{\cT_i}(X,Z)
$$
were all specializations from $A$. This holds
even true if we consider any finite family of of pairs $(X,Y)$
simultaneously.

In the following examples, we take for $\Delta$ the universal degree
function with $A=K(\cA)$ and put $\cT(\cA):=\cT(\cA,\Delta)$. We
choose a field $K$ of characteristic zero.

\Examples: \Bsp. [Deligne's case] Let $\cA=\cB^{\|op|}$ where
$\cB=\FS$ is the category of finite sets. Then $K(\cA)=\ZZ[t]$. The
uniform functors are of the form $P(A):=\|Hom|_\FS(A,X)$ where $X$ be
a finite set. The adapted degree function corresponds to
$t=n:=|X|$. Since $\|Aut|_\FS(X)\cong S_n$ we see that {\it $\cT(\cA)$
interpolates the categories $\|Rep|(S_n,K), n\ge0$}.

\Bsp. More generally, let $G$ be a finite group and let
$\cA=\cB^{\|op|}$ where $\cB$ is the category of finite
sets with {\it free} $G$\_action. Then $K(\cA)=\ZZ[t]$. For any object
$X$ of $\cA$ we get a uniform functor $P(A)=\|Hom|_\cB(A,X)$
whose adapted degree function has parameter $t=|X|=|G|n$, where
$n=|X/G|$. We have $\|Aut|_\cA(X)\cong S_n\wr G=S_n\semidir G^n$. Thus
{\it $\cT(\cA)$ interpolates $\|Rep|(S_n\wr G,K)$, $n\ge0$.} In
particular, for $G=\ZZ_2$, we get the representation categories of the
hyperoctahedral groups, i.e., the Weyl groups of type ${\Ss BC}_n$.

\Bsp. Let $\cA=\cB^{\|op|}$ where $\cB$ is the category of chains
$$
A_1\leftarrow A_2\leftarrow A_3\leftarrow A_4\leftarrow\ldots
$$
of finite sets with $A_i=\emptyset$ for $i>\!\!>0$. This category is
also equivalent to the category of finite rooted trees with graph maps
which preserve the distance to the root. Then
$K(\cA)=\ZZ[t_1,t_2,\ldots]$. The uniform pro\_objects are chains of
finite sets
$$
X_1\leftarrow X_2\leftarrow X_3\leftarrow X_4\leftarrow\ldots
$$
(possibly non\_empty for all $i$) such that all connecting maps are
uniform. The corresponding parameters are $t_i=|X_i|=:n_i$. Thus {\it
$\cT(\cA)$ interpolates $\|Rep|(S_{n_1}\wr (S_{n_2}\wr
(S_{n_3}\wr S_{n_4}\wr\ldots))),K)$, with $n_i\ge1$.}

\Bsp. Let $\cA$ be the category of finite solvable groups. Then
$K(\cA)=\ZZ[t_p\mid p\ {\rm prime}]$. Let $FS_n$ be the free
pro-solvable group on $n$ letters (i.e., the completion of the free
group $F_n$ on $n$ letters with respect to the topology defined by all
normal subgroups $N$ such that $F_n/N$ is solvable). This is a uniform
pro\_object. Its adapted degree function has $t_p=p^n$. A moments
thought shows that these functions are Zariski\_dense in
$\|Spec|K(\cA)$. Thus {\it $\cT(\cA)$ interpolates
$\|Rep|(\|Aut|FS_n,K)$, $n\ge0$.}

\Bsp. Let $\cA=\FV_{\FF_q}$, the category of finite dimensional
$\FF_q$\_vector spaces. Then $K(\cA)=\ZZ[t]$. Every object $X$ of
$\cA$ is uniform with adapted degree function corresponding to the
parameter $t=|X|=q^n$. Thus {\it $\cT(\cA)$ interpolates
$\|Rep|(GL_n(\FF_q),K)$, $n\ge 0$, $q$ fixed.}

\Bsp. Let $\cA$ be the category of chains of homomorphisms of finite
dimensional $\FF_q$\_vector spaces
$$
\ldots\pfeil V_{-2}\pfeil V_{-1}\pfeil V_0\pfeil V_2\pfeil
V_3\pfeil V_4\pfeil\ldots
$$ with $V_i=0$ for $|i|>\!\!>0$. Then
$K(\cA)=\ZZ[\ldots,t_{-1},t_0,t_1,\ldots]$. Moreover, a uniform
pro\_object is also a chain but all homomorphisms have to be injective
and the vanishing condition on the $V_i$ has to be replaced by the
condition that the inverse limit is zero. In other words, a uniform
object is an $\FF_q$\_vector space $X$ with a $\ZZ$\_filtration $X_i$
such that $X=\bigcup_iX_i$ and $0=\bigcap_iX_i$. This shows that
{\it $\cA(\cA)$ interpolates
$\|Rep|(P_{\ldots n_{-1}n_0n_1\ldots},K)$.} Here $P_{\ldots
n_{-1}n_0n_1\ldots}$ is group of invertible block matrices
$$
\pmatrix{\ddots\cr
&B_{-1-1}&B_{-10}&B_{-11}\cr
&0&B_{00}&B_{01}\cr
&0&0&B_{11}\cr
&&&&\ddots\cr}
$$
where $B_{ij}$ is an $n_i\times n_j$\_matrix with entries in $\FF_q$.

\Bsp. Let $\cA$ be the category of (non\_empty) affine spaces. Then
$K(\cA)=\ZZ[t]$ and every object is uniform. Thus {\it $\cT(\cA)$
interpolates $\|Rep|(A_n,K)$, $n\ge0$} where $A_n=GL(n,\FF^q)\semidir\FF_q^n$
is the affine group.

\Appendix A. Protomodular and Mal'cev categories

\advance\appcount1\SATZ1\GNo0 \edef\Aname{\AppName\number\appcount}%
In this appendix, we recall two classes of categories, namely
protomodular and Mal'cev ones, which
generalize various aspects of the category of groups and abelian
categories. Mal'cev categories are more general than protomodular
ones.

In a nutshell, the short five\_lemma holds for protomodular categories
while some version of Jordan\_H\"older theorem can be proved for
Mal'cev categories (see Appendix~B). Moreover, Mal'cev categories have
a nice theory of relations. That's why they are of particular interest
to us. It should be pointed out, though, that most, if not all,
natural examples are already protomodular.

\Definition Malcevdef. A regular category is a {\it Mal'cev
category} if every reflexive relation is an equivalence relation.

\noindent This definition is due to Carboni, Lambek, and
Pedicchio~\cite{CLP}. It is the following property of exact Mal'cev
categories which we are actually going to use.

\Proposition pushpull. {\rm(\cite{CKP} Thm.~5.7)} Let $\cA$ be an
exact Mal'cev category. Then every pair of surjective morphisms with
the same domain has a push-out. Moreover for any commutative diagram
$$2
\cxymatrix{
u\ar@{>>}[d]\ar@{>>}[r]&y\ar@{>>}[d]\\
x\ar@{>>}[r]&z
}
$$
of surjective morphisms the following are equivalent:
\Item{27}Diagram \cite{E2} is a pull\_back.
\Item{28}Diagram \cite{E2} is a push\_out and $u\pfeil x\times y$ is
injective.\Par

\Remark: The implication \cite{I27}$\Rightarrow$\cite{I28} is a
general property of regular categories. It is the characterization of
pull\_backs in term of push\_outs which is of particular interest.

Now, we define protomodular categories. They will be used in
\cite{genericdegree} stating the non\_singularity of the generic
degree function.  First, recall the category $\|Pt|_s\cA$ of triples
$(x,e,d)$ where $e:x\pfeil s$ and $d:s\pfeil x$ are morphisms with
$ed=1_s$. It is a pointed category with zero object $(s,1_s,1_s)$.

\Definition protodef. A regular category $\cA$ is {\it protomodular}
if for any morphism $\sq\pfeil s$ the pull\_back functor
$\|Pt|_s\cA\pfeil\|Pt|_\sq\cA$ reflects isomorphisms.

\noindent In protomodular categories we have the following form of the
short five\_lemma:

\Lemma fivelemma. Let $\cA$ be a pointed protomodular
category and consider the following commutative diagram:
$$
\cxymatrix{
x\inj[r]^i\ar[d]^a&y\sur[r]^e\ar[d]^b&z\ar[d]^c\\
x'\inj[r]^{i'}&y'\sur[r]^{e'}&z'
}
$$
where $i,i'$ is the kernel of $e,e'$, respectively. Assume that $a$
and $c$ are isomorphisms. Then $b$ is an isomorphism, as well.

\excount=0 \Remark: \Bsp. For further reading we recommend the book
\cite{BoBo} by Borceux and Bourn.

\Bsp. Our notion of a Mal'cev/protomodular category
differs slightly from the one in the literature (like, e.g.,
\cite{BoBo}) since these categories are usually not required to be
regular. Instead, a certain amount of limits is required to exist. For
us, it doesn't really matter since we are only interested in regular
categories anyway. A sufficient amount of limits is already build
into them.

\Examples: \Bsp. Every abelian category is protomodular since in that case
the pull\_back functor $\|Pt|_s\cA\pfeil\|Pt|_\sq\cA$ is even an
equivalence of categories.

\Bsp. Let $\cA$ be the category of models of an equational
theory. Then $\cA$ is always exact. If the theory contains a group
operation then $\cA$ is protomodular. Therefore, the categories of
groups, rings (with or without unity), Boolean algebras, Lie algebras,
or any other type of algebras is protomodular. The full story is as
follows (see \cite{BoBo} Thm.~3.1.6): $\cA$ is protomodular if and
only if there is an $n\in\NN$ such that the theory contains $n$
constants $e_1,\ldots,e_n$, $n$ binary operations
$a_1(x,y),\ldots,a_n(x,y)$, and one $(n+1)$\_ary operation
$b(x_0,\ldots,x_n)$ such that
$$
a_1(x,x)=e_1,\ldots,a_n(x,x)=e_n,\ b(x,a_1(x,y),\ldots,a_n(x,y))=y.
$$
For example, for groups one has $n=1$, $e_1=e$, $a_1(x,y)=x^{-1}y$ and
$b(x,y)=xy$.

\Bsp. The opposite of the category of sets or, more generally, of an
elementary topos is exact protomodular (\cite{BoBo} Ex.~3.1.17).

\Bsp. Every protomodular category is Mal'cev (see, e.g.,
\cite{BoBo} Prop.~3.1.19).

\Bsp. Let $\cA$ be the category of models of an equational
theory. Then $\cA$ is Mal'cev if and only if the theory contains a
ternary operation $m(x,y,z)$ with $m(x,x,z)=z$ and $m(x,z,z)=x$ for
all $x,z$ (\cite{BoBo} Thm.~2.2.2). Any group structure gives rise to
such an operation, namely $m(x,y,z)=xy^{-1}z$. If $\cA$ is the
category of sets equipped with such a ternary operation then $\cA$ is
Mal'cev but not protomodular. This is shown in
\S\cite{determinant}~example~9.

\Bsp. The categories of (finite) sets, (finite) monoids, (finite)
posets, and (finite) lattices are not Mal'cev.

Let $\fP$ be one of the properties ``protomodular'' or
``Mal'cev''. Then the class of $\fP$\_categories enjoys many
permanence properties. In the following, let $\cA$ be any
$\fP$\_category.

\Bsp. Any full subcategory of $\cA$ which is closed under products,
subobjects and quotients is again $\fP$. Thus the category of finite
models for an equational $\fP$\_theory is again $\fP$. Examples are
the categories of finite groups, finite rings, and finite Boolean
algebras. The latter example is, by the way, equivalent to $\FS^{\rm
op}$, the opposite category of finite sets.

\Bsp. Let $\cD$ be a small category. Then the diagram category
$[\cD^{\rm op},\cA]$ is $\fP$. This applies, e.g., to the category of
arrows $x\pfeil y$ in $\cA$ or the category of objects with
$G$\_action where $G$ a fixed group (or monoid).

\Bsp. Fix an object $s$ of $\cA$. Then the slice category $\cA/s$ of
$s$\_objects $x\pfeil s$ is $\fP$. The same holds for the coslice category
$s\backslash\cA$ of arrows $s\pfeil x$ and the category of points
$\|Pt|_s\cA=(s\backslash\cA)/(s\pfeil s)$.

\Bsp. The category $\cA\mod s$ of ``dominant'' $s$\_objects $x\auf s$
is $\fP$. The same holds for the full subcategory category
$s\pfeil\cA$ of objects such that there exists an arrow $s\pfeil x$
provided $s$ is projective in it.

\Appendix B. Degree functions on Mal'cev categories

In this appendix, we determine the degree functions on certain Mal'cev
categories. We also state a result to the effect that degree functions
separate certain morphisms. For simplicity, we restrict to
categories~$\cA$ which are essentially small. This has the effect that
$\cA$ has a universal degree function $\Delta:\fE(\cA)\pfeil K(\cA)$
(see \cite{E81}).

\appsection Degreefunctions. Results

In this section, we are only stating the results. Proofs are given in
section~B3.

First, we need the following finiteness condition.

\Definition. A regular category is of {\it finite type} if $\|sub|(x)$
satisfies for every object $x$ the ascending and the descending chain
condition.

\noindent This condition implies, in particular, that also the set of
quotient objects of any $x$ satisfies the ascending and descending chain
condition.

In determining degree functions, we first consider pointed
categories. This means that there is an object $\0$ which is both
initial and terminal. An object $x\not\cong\0$ is {\it simple} if, up
to isomorphism, $\0$ and $x$ are its only quotients.

\Theorem degreepointed. Let $\cA$ be an essentially small, pointed, exact
Mal'cev category of finite type. Let $S(\cA)$ be its set of
isomorphism classes of simple objects. Then the map
$s\mapsto\<s\auf\0\>$ induces an isomorphism $\ZZ[S(\cA)]\pf\sim
K(\cA)$.

The case of non\_pointed categories can be often reduced to the
pointed case. Let $s$ be an object of $\cA$. Recall that
$\|Pt|_s\cA$ is the category of triples $(x,e,d)$ where $e:x\pfeil
s$ and $d:s\pfeil x$ are morphisms with $ed=1_s$. This is a pointed
category with zero object $(s,1_s,1_s)$.

\Theorem degreeinitial. Let $\cA$ be an essentially small, exact
Mal'cev category of finite type. Assume that $\cA$ has an initial
object $\0$. Then the forgetful functor $\|Pt|_\0\cA\pfeil\cA$ induces
an isomorphism $K(\|Pt|_\0\cA)\pf\sim K(\cA)$.\footnote{Observe
$\|Pt|_\0\cA=\cA/\0$ since $s:\0\pfeil x$ is redundant.}

If $\cA$ doesn't have an initial element then we look at {\it minimal}
objects, i.e., objects which have no proper subobjects. Equivalently,
an object $m$ is minimal if every morphism $x\pfeil m$ is
surjective. If $\cA$ is of finite type then every object has a minimal
subobject. If, moreover, $\cA$ has all finite limits (hence
intersections) then this minimal subobject $x_{\|min|}$ is even
unique. Another consequence is that for any two minimal objects $m$
and $\mq$ there is at most one morphism $\mq\pfeil m$. In fact, the
graph of this morphism would have to be $(\mq\times m)_{\|min|}$. The
set $M(\cA)$ of isomorphism classes of minimal objects is a
join-semilattice with $\mq\ge m$ if there is a morphism $\mq\pfeil m$
and $m\vee\mq=(m\times\mq)_{\|min|}$. If there is a morphism
$\mq\pfeil m$ then there is a pull\_back functor $\Phi_{\mq\pfeil
m}:\|Pt|_m\cA\pfeil\|Pt|_\mq\cA$.

\Theorem degreecomplete. Let $\cA$ be a essentially small, complete,
exact Mal'cev category of finite type.
\Item{74}Let $m,\mq$ be
two minimal objects. Then $\Phi_{\mq\pfeil m}$ preserves simple
objects. In particular we can define
$$
S(\cA):=\|lim|_{m\in M(\cA)}S(\|Pt|_m\cA).
$$
\Item{75}The map $(s\lrarrow m)\mapsto\<s\auf m\>$ induces an
isomorphism $\ZZ[S(\cA)]\pf\sim K(\cA)$.\Par

\Example: Let $\cA$ be the category of finite unital commutative
rings. Then a ring $A$ is minimal if $A=\ZZ_n$ is cyclic. Via the
correspondence, $I\mapsto I\oplus\ZZ_n$, the category
$\|Pt|_{\ZZ_n}\cA$ is equivalent to the category of finite
non\_unital rings $I$ with $nI=0$. The pull\_back functor $\Phi$ sends
$I\oplus\ZZ_n$ to $I\oplus\ZZ_m$ (provided $n|m$). Thus, $S(\cA)$ is
the set of isomorphism classes of finite simple non\_unital
commutative rings. They are easily classified: let $I$ be finite
simple commutative. Because of simplicity we have $xI=0$ or $xI=I$ for
all $x\in I$. For the same reason, the set $J=\{x\in I\mid xI=0\}$ is
either $0$ or $I$. If $J=I$ then $II=0$. Thus $I=N_p$ where $p$ is a
prime and $N_p=\ZZ_p$ as an additive group but with zero
multiplication. If $J=0$ then multiplication by any $x\ne0$ is
surjective. Thus there is $e\in I$ with $xe=x$. But then $ze=z$ for
all $z\in xI=I$, i.e., $e\in I$ is an identity element. This implies
that $I=\FF_q$ is a finite field. Thus $K(\cA)$ is a polynomial ring
over two sets of variables $n_p$, $p$ prime and $f_q$, $q$ a prime
power.

\medskip

An application of the theory above is the following statement on the
values of the universal degree functions.

\Theorem protoDelta. Let $\cA$ be an essentially small, complete,
exact protomodular category of finite type with universal degree
function $\Delta:\fE(\cA)\pfeil K(\cA)$. Consider the commutative
diagram
$$83
\cxymatrix{u\inj[r]^j\sur[rd]_\eq&x\sur[d]^e\\&y}.
$$
Then $\Delta(e)=\Delta(\eq)$ if and only if $j$ is an isomorphism.

\Remark: Example~9 at the end of \S\cite{determinant} shows that the
theorem fails for Mal'cev categories.

\appsection JordanHolder. Lambek's Jordan-H\"older theorem

In this section (only), {\it the classical product of relations will
be denoted by $rs$ instead of $r\circ s$.}

Let $\cA$ be a {\it pointed} regular category. Then we can talk about
the {\it kernel} $\|ker|f=x\times_y\0$ of a morphism $f:x\pfeil y$. A
{\it normal series} of an object $x$ is a diagram
$$65
\cxymatrix{
\0\inj[r]&x_1\inj[r]\sur[d]^\sim&x_2\inj[r]\sur[d]&
\cdots\inj[r]&x_{n-1}\inj[r]\sur[d]&
x_n\sur[d]\ar@{}[r]|{=}&x\\
&y_1&y_2&&y_{n-1}&y_n}
$$
where each (horizontal) injective morphism is the kernel of the
following (vertical) surjective morphism. The $y_i$ are called the
{\it factors} of the normal series. Two series are {\it equivalent}
if, after a suitable permutation, their factors are isomorphic.

Given a normal series of one of the factors $y_i$ one can {\it refine}
the normal series of $x$:
$$87
\cxymatrix{
\cdots\inj[r]&x_{i-1}\inj[r]\sur[d]&x_1'\inj[r]\sur[d]\ar@{.>>}@/^1pc/[dd]&
x_2'\inj[r]\sur[d]\ar@{.>>}@/^1pc/[dd]&\cdots\inj[r]&
x_{m-1}'\inj[r]\sur[d]\ar@{.>>}@/^1pc/[dd]&x_i\inj[r]\sur[d]\ar@{.>>}@/^1pc/[dd]&
x_{i+1}\inj[r]\sur[d]&\cdots\\
&y_{i-1}&y_1'\inj[r]\sur[d]^\sim&y_2'\inj[r]\sur[d]&
\cdots\inj[r]&y_{m-1}'\inj[r]\sur[d]&y_i\sur[d]&y_{i+1}\\
&&z_1&z_2&&z_{m-1}&z_m}
$$
Here all squares are pull\_backs. Of course one can refine all $y_i$
simultaneously giving the most general refinement. Now we have the
following Schreier type theorem:

\Theorem Schreier. Let $x$ be an object of a pointed, exact Mal'cev
category $\cA$. Then any two normal series of $x$ have equivalent
refinements.

The theorem is essentially due to Lambek, \cite{La1}, \cite{La2},
except that he proves it just for equational theories and that our
notion of refinement seems to be stricter than his. Therefore, we
repeat his proof and observe that it stays valid.

First, we reformulate the theorem in terms of relations and state a
series of lemmas. Every subquotient diagram
$\xymatrix@-10pt{x&y\inj[l]\sur[r]&z}$ gives rise to an equivalence
relation $r$ on $y$ which we may consider as a subobject of $x\times
x$. Conversely, the subquotient is uniquely determined by $r$: for any
subobject $u$ of $x$ we denote the image of the first projection
$r\times_xu\pfeil x$ by $ru$. Then $y=rx$ and $z=y/r$. Since $\cA$ is
exact, the relations which arise as $r$ are precisely the {\it
subequivalence relations}, i.e., the relations which are symmetric and
transitive but not necessarily reflexive.

A normal series of $x$ can be encoded as a sequence
$r_1,r_2,\ldots,r_n$ of subequivalence relations where $r_i$
corresponds to the subquotient diagram $\xymatrix@=10pt{x&x_i\ar@{
>->}[l]\ar@{>>}[r]&y_i}$. In other words, $r_i=x_i\times_{y_i}x_i\into
x\times x$. Since $x_i=r_ix$ and $\|ker|(x_i\auf
y_i)=x_i\0$, the kernel conditions are equivalent to
$$
\0=r_1\0,r_1x=r_2\0,\ldots,r_{n-1}x=r_n\0,r_nx=x.
$$
Next, we need to encode refinements:

\Definition. A {\it refinement} of a subequivalence
relation $r$ on $x$ is a sequence of subequivalence
relations $s_1,\ldots,s_m$ with
$$0
\eqalignno{%
61&&rs_jr=s_j,\quad j=1,\ldots,m\quad{\rm and}\cr
66&&r\0=s_1\0,\ s_1x=s_2\0,\ s_2x=s_3\0,\ldots,
s_{m-1}x=s_m\0,\ s_mx=rx.\cr
}
$$

\Remark: It is condition \cite{E61} which is missing in Lambek's
notion of refinement.

\Lemma. There is a natural one-to-one correspondence between
refinements of $r$ and normal series of $rx/r$.

\Proof: Consider the quotient $rx\auf r/x$ and let $u$ be a subobject
of $x$. Then it is well known that $u$ is the pull\_back of a
subobject of $rx/r$ if and only if $ru=u$. If we apply this to the
morphism $x\times x\auf x/r\times x/r$ then we obtain that a relation
$s$ on $x$ is the pull\_back of a relation on $x/r$ if and only if
$rsr=s$. Thus, \cite{E61} makes asserts that the $s_j$ are pull\_backs
from subequivalence relations $s_j'$ on $y:=rx/r$. The conditions
\cite{E66} imply that the $s_j'$ form a normal series of $y$.\qed

\Lemma MalcevLemma1. Let $\cA$ be a pointed Mal'cev category and $r$,
$s$ two subequivalence relations on $x$. Then $rsr\0=rs\0$ and
$rsrx=rsx$.

\Proof: Since $\cA$ is Mal'cev, for any relation $p\subseteq x\times
y$ holds
$$
pp^\vee p=p.
$$
Applying this to $p=rs$ we get
$$88
rsrs=rs.
$$
Thus,
$$0
\eqalignno{
&&rs\0\subseteq rsr\0\subseteq rsrs\0=rs\0\cr
&&rsx=rsrsx\subseteq rsrx\subseteq rs x.\cr}
$$\qed

\noindent{\it Proof of \cite{Schreier}:} Let $r_1,\ldots,r_n$ and
$s_1,\ldots,s_m$ encode two normal series and put $r_{ij}:=r_is_jr_i$,
$i=1,\ldots,n$, $j=1,\ldots,m$. For a fixed $i$ we claim that
$r_{i1},\ldots,r_{im}$ refines $r_i$. First, $r_{ij}$ is clearly
symmetric. Transitivity holds by \cite{E88}:
$$
r_{ij}r_{ij}=r_is_jr_is_jr_i=r_is_jr_i=r_{ij}. 
$$
Hence $r_{ij}$ is a subequivalence relation. Condition \cite{E61}
holds trivially. Finally, we have, using \cite{MalcevLemma1},
$$0
\eqalignno{
&&r_{i1}\0=r_is_1r_i\0=r_is_1\0=r_i\0\cr
&&r_{ij}\0=r_is_jr_i\0=r_is_j\0=r_is_{j-1}x=r_is_{j-1}r_ix=r_{ij-1}x,\quad
j=2,\ldots,m\cr
&&r_{im}x=r_is_mr_ix=r_is_mx=r_ix.\cr
}
$$

Symmetrically, put $s_{ji}:=s_jr_is_j$. Then $s_{j1},\ldots,s_{jn}$
forms a refinement of $s_j$. The assertion follows now from Lambek's
version of the Zassenhaus butterfly lemma (\cite{La1} Prop.~3, Thm.~I,
\cite{CLP} Prop.~4.2):
$$
r_{ij}x/r_{ij}\cong s_{ji}x/s_{ji}.
$$\qed

A normal series is called a {\it composition series} if all of its
factors are simple. Now we have the following Jordan\_H\"older type
theorem:

\Theorem JH. Let $\cA$ be a pointed, exact Mal'cev category of finite
type. Then every object has a composition series and any two
composition series are equivalent.

\Proof: The additional condition \cite{E61} makes sure that a
composition series can be only refined in a trivial way. In fact, if
$y_i$ is simple then the only possible refinements are
$$
\cxymatrix{
\0\inj[r]\sur[d]&\cdots\inj[r]&\0\inj[r]\sur[d]&
y_i\inj[r]^\sim\sur[d]^\sim&y_i\inj[r]^\sim\sur[d]&
\cdots\inj[r]^\sim&y_i\sur[d]\\
\0&\cdots&\0&y_i&\0&\cdots&\0}
$$
Thus all non\_zero factors stay the same. Now the assertion follows
from \cite{Schreier}.\qed

\appsection Proofs. Proofs of Theorems \ncite{degreepointed}\kern2pt--\ncite{protoDelta}

Let $\cA$ be a pointed regular category and $\fS(\cA)$ its class of
simple objects. A {\it rank function} on $\cA$ is a function
$\rho:\|Ob|\cA\pfeil K$ such that
$$
\rho(x)=\rho(\|ker|e)\cdot\rho(y)\quad\hbox{for all surjective $e:x\auf y$.}
$$

\Proposition extend4. Let $\cA$ be a pointed, exact Mal'cev category of
finite type. Then every function $\rho_0:S(\cA)\pfeil K$
extends uniquely to a rank function $\rho:\|Ob|\cA\pfeil K$.

\Proof: For any object $x$ we choose a composition series with
simple factors $y_1,\ldots, y_n$. Then we are forced to define
$\rho(x)=\prod_{i=1}^n\rho_0(y_i)$, showing the uniqueness of the
extension. Moreover, \cite{JH} shows that $\rho(x)$ is well
defined. If $e:x\auf y$ is surjective then the composition factors of
$x$ are those of $y$ together with those of $\|ker|e$ which shows that
$\rho$ is a rank function.\qed

\noindent{\it Proof of \cite{degreepointed}:} In view of \cite{extend4}
it suffices to show that
$$64
\delta(e):=\rho(\|ker|e),\qquad\rho(x):=\delta(x\auf\0)
$$
establishes a bijection between rank and degree functions on $\cA$.

First we start with $\delta$. Then $\rho$ as in \cite{E64} is a rank
function. Indeed, let $f:x\auf y$ be surjective. Then
$$
\eqalign{\rho(x)&=\delta(x\auf y\auf\0)=\delta(x\auf
y)\delta(y\auf\0)=\cr
&=\delta(\|ker|x\auf\0)\delta(y\auf\0)=\rho(\|ker|f)\rho(y).\cr}
$$

Conversely, given $\rho$ then $\delta$ is a degree function. Indeed,
the invariance under pull\_backs is clear. Multiplicativity under
composition follows from the following diagram where all squares are
pull\_backs.
$$
\cxymatrix{
\|ker|e'\inj[d]\sur[r]&\0\inj[d]\\
\|ker|ee'\inj[d]\sur[r]&\|ker|e\inj[d]\sur[r]&\0\inj[d]\\
x\sur[r]^{e'}&y\sur[r]^e&z
}
$$
It is clear that both assignments are inverse to each other.\qed

\noindent{\it Proof of \cite{degreeinitial}:} First observe that since
$\0$ is initial, the arrow $\0\pfeil x$ is redundant and the category
$\|Pt|_\0\cA$ is the same as the slice category $\cA/\0$. Now we
can define maps between $K(\cA)$ and $K(\|Pt|_0\cA)$
$$
\eqalign{
K(\cA)\pfeil K(\|Pt|_\0\cA):&\ \ \<x\auf y\>\mapsto\<\xymatrix@=4pt{x\times_y\0\sur[rr]\ar[rd]&&\0\ar[ld]\\&\0}\>\cr
K(\|Pt|_\0\cA)\pfeil K(\cA):&\ \
\<\xymatrix@=4pt{x\sur[rr]\ar[rd]&&y\ar[ld]\\&\0}\>\mapsto\<x\auf y\>\cr
}
$$
which are clearly inverse to each other.\qed

\noindent{\it Proof of \cite{degreecomplete}:} Consider the following
pull\_back diagram
$$79
\cxymatrix{
\sq\ar@{>>}[r]^\eq\ar@{>>}[d]&\mq\ar@{>>}[d]\\
s\ar@{>>}[r]^e&m
}
$$
where $s$ is simple in $\|Pt|_m\cA$, i.e., $e$ is indecomposable in
$\cA$. We have to show that $\eq$ is indecomposable, as well. The
morphism $\eq$ is not an isomorphism since diagram \cite{E79} is also
a push\_out diagram (see remark after \cite{pushpull}). Now assume
that $\eq$ factorizes. Then we get the diagram
$$78
\cxymatrix{
\sq\ar@{>>}[r]\ar@{>>}[d]&u\ar@{>>}[r]&\mq\ar@{>>}[d]\\
s\ar@{>>}[rr]&&m
}
$$
Since $\cA$ is an exact Mal'cev category, the quotient objects of $x$
form a {\it modular} lattice (\cite{CKP} Prop.~3.3). We have
$s\vee\mq=\sq$ since \cite{E79} is a pull\_back. The indecomposability
of $s\auf m$ implies either $u\wedge s=s$ or $u\wedge s=m$. Thus
$$
u=u\wedge\sq=u\wedge(s\vee\mq)=(u\wedge s)\vee\mq=
\cases{s\vee\mq=\sq&or\cr m\vee\mq=\mq.&\cr} 
$$
Thus $\sq$ is simple as well.

\noindent\cite{I75}: If $m$ is a minimal object of $\cA$ let
$m\backslash\cA$ be the coslice category of all arrows $m\pfeil
x$. Since any morphism $m\pfeil x$ is necessarily unique this is a
full subcategory of $\cA$ with initial element $m$. It is clear that
$m\backslash\cA$ is again an exact Mal'cev category. Moreover, $\cA$
is the union of all subcategories $m\backslash\cA$ where $m$ runs
through all minimal objects.

For a morphism $\mq\auf m$ of minimal objects consider the commutative
diagram
$$80
\cxymatrix{
K(\|Pt|_m\cA)\ar[d]_{\Phi_{\mq\auf m}}\ar[r]^\sim&K(m\backslash\cA)\ar[d]\\
K(\|Pt|_\mq\cA)\ar[r]^\sim&K(\mq\backslash\cA)
}
$$
The two horizontal arrows are isomorphisms by
\cite{degreepointed}. The right vertical arrow comes from the
inclusion $m\backslash\cA\subseteq\mq\backslash\cA$. It is easy to
check that the diagram commutes. Thus, we have to show that
$$82
\|lim|_{m\in M(\cA)}K(m\backslash\cA)\Pf\sim K(\cA).
$$
Let $e:x\auf y$ be surjective. Then $e$ is in the image of
$K(m\backslash\cA)$ with $m=x_{\|min|}$ showing that \cite{E80} is
surjective. On the other hand, let $e_i\in\fE(m_i\backslash\cA)$ with
$\<e_1\>=\<e_2\>$ in $K(\cA)$. To show this inequality only finitely
many objects and morphisms of $\cA$ are needed. Thus, this equality
holds already in $K(m\backslash\cA)$ with $m\in M(\cA)$ big
enough. This shows that \cite{E80} is injective.\qed

\noindent{\it Proof of \cite{protoDelta}}: Pull\_back by the
surjective morphism $\eq$ reflects isomorphisms and preserves
$\Delta$. Thus we may assume that $\eq$ has a splitting $s$. Thus the
diagram takes place in $\|Pt|_y\cA$. By \cite{E82} there is a minimal
object $m\ge y_{\|min|}$ such that $e$ and $\eq$ have the same image
in $m\backslash\cA$. This means that also the pull\_backs of $e$,
$\eq$ by $m\pfeil y$ have the same image in $\|Pt|_m\cA$. On the 
hand, since $\cA$ is protomodular, the pull\_back functor
$\|Pt|_y\cA\pfeil\|Pt|_m\cA$ reflects isomorphisms. Thus it suffices
to prove the assertion for $\|Pt|_m\cA$, i.e., we may assume from now
on that $\cA$ is pointed and $y=0$.

Consider a composition series \cite{E65} of $x$. Then the
$u_i:=x_i\cap u$ form a normal series of $u$ with factors
$z_i:=\|image|(x_i\cap u\pfeil y_i)$. I claim that $z_i=y_i$ for all
$i$. For that, we construct a directed graph. The vertices are the $i$
with $z_i\ne y_i$. We draw an arrow $i\pfeil j$ if $y_i$ is isomorphic
to a composition factor of $z_i$. In that case
$$84
y_i\preceq z_j\preceq y_j.
$$
The assumption $\Delta(u\auf\0)=\Delta(x\auf\0)$ means that $u$ and
$x$ have the same composition factors. This implies that each vertex
has at least one outgoing edge. Therefore, the graph must contain a
directed cycle. If $j$ is part of such a cycle then
\cite{partialorder0} and \cite{E84} imply $z_j=y_j$,
contradicting the choice of $j$. This proves the claim.

Finally, we show by induction on $n$, the number of composition
factors, that $u=x$. For that consider the following diagram
$$
\cxymatrix{
u_{n-1}\inj[r]\ar[d]^\sim&u_n\sur[r]\ar[d]&z_n\ar[d]^\sim\\
x_{n-1}\inj[r]&x_n\sur[r]&y_n
}
$$
Here the left horizontal arrows are the kernels of the right
horizontal ones. The left vertical arrow is an isomorphism by
induction. The right vertical arrow is an isomorphism by what we
showed above. \cite{fivelemma} implies $u=u_n\pf\sim x_n=x$.\qed

\rightskip0pt plus 5pt minus 5pt

\beginrefs

\L|Abk:AKO|Sig:AKO|Au:Andr\'e, Y.; Kahn, B.; O'Sullivan, P.|Tit:Nilpotence,
radicaux et structures mono\"\i\-dales|Zs:Rend. Sem. Mat. Univ.
Padova|Bd:108|S:107--291|J:2002|xxx:math.CT/0203273||

\B|Abk:AM|Sig:AM|Au:Artin, M.; Mazur, B.|Tit:Etale
homotopy|Reihe:Lecture Notes in Mathematics
{\bf100}|Verlag:Springer-Verlag|Ort:Berlin-New York|J:1969|xxx:-||

\B|Abk:Barr|Sig:Ba|Au:Barr, M.|Tit:Exact categories|Reihe:Lecture
Notes in Mathematics
{\bf236}|Verlag:Springer-Verlag|Ort:Berlin|J:1971|xxx:-||

\B|Abk:Borceux|Sig:Bo|Au:Borceux, F.|Tit:Handbook of categorical
algebra. 2. Categories and structures|Reihe:Encyclopedia of
Mathematics and its Applications, {\bf51}|Verlag:Cambridge University
Press|Ort:Cambridge|J:1994|xxx:-||

\B|Abk:BoBo|Sig:BB|Au:Borceux, F.; Bourn, D.|Tit:Mal'cev,
protomodular, homological and semi\_abelian
categories|Reihe:Mathematics and its Applications {\bf
566}|Verlag:Kluwer Academic Publishers|Ort:Dordrecht|J:2004|xxx:-||

\L|Abk:CKP|Sig:CKP|Au:Carboni, A.; Kelly, G.; Pedicchio, M.|Tit:Some
remarks on Maltsev and Goursat categories|Zs:Appl. Categ.
Structures|Bd:1|S:385--421|J:1993|xxx:-||

\L|Abk:CLP|Sig:CLP|Au:Carboni, A.; Lambek, J.; Pedicchio,
M.|Tit:Diagram chasing in Mal'cev categories|Zs:J. Pure Appl.
Algebra|Bd:69|S:271--284|J:1991|xxx:-||

\L|Abk:De|Sig:De|Au:Deligne, P.|Tit:La cat\'egorie des
repr\'esentations du groupe sym\'etrique $S_t$ lorsque $t$ n'est pas
un entier naturel|Zs:Preprint|Bd:-|S:78
pages|J:-|xxx:www.math.ias.edu/\lower4pt\hbox{\char126}phares/deligne/Symetrique.pdf||

\Pr|Abk:EKMS|Sig:EKMS|Au:Ern\'e, M.; Koslowski, J.; Melton, A.;
Strecker, G.|Artikel:A primer on Galois connections|Titel:Papers on
general topology and applications (Madison, WI, 1991)|Hgr:-|Reihe:Ann.
New York Acad. Sci.|Bd:704|Verlag:New York Acad. Sci.|Ort:New
York|S:1993|J:103--125|xxx:www.iti.cs.tu-bs.de/TI-INFO/koslowj/RESEARCH/gal\string_bw.ps.gz||

\L|Abk:Gou|Sig:Gou|Au:Goursat, E.|Tit:Sur les substitutions
orthogonales et les divisions r\'eguli\`eres de
l'espace|Zs:Ann. Sci. \'Ecole
Norm. Sup. (3)|Bd:6|S:9--102|J:1889|xxx:-||

\L|Abk:Greene|Sig:Gr|Au:Greene, C.|Tit:On the M\"obius algebra of a
partially ordered set|Zs:Advances in
Math.|Bd:10|S:177--187|J:1973|xxx:-||

\L|Abk:KnTen|Sig:Kn|Au:Knop, F.|Tit:A construction of semisimple
tensor categories|Zs:C. R. Math. Acad. Sci.
Paris|Bd:343|S:15--18|J:2006|xxx:math.CT/0605126||

\L|Abk:La1|Sig:\\La|Au:Lambek, J.|Tit:Goursat's theorem and the
Zassenhaus lemma|Zs:Canad. J. Math.|Bd:10|S:45--56|J:1958|xxx:-||

\Pr|Abk:La2|Sig:\\La|Au:Lambek, J.|Artikel:On the ubiquity of
Mal'cev operations|Titel:Proceedings of the International
Conference on Algebra, Part 3 (Novosibirsk,
1989)|Hgr:-|Reihe:Contemp. Math.|Bd:131, \rm Part
3|Verlag:Amer. Math. Soc.|Ort:Providence, RI|S:135--146|J:1992|xxx:-||

\L|Abk:Li|Sig:Li|Au:Lindstr\"om, B.|Tit:Determinants on
semilattices|Zs:Proc. Amer Math. Soc.|Bd:20|S:207--208|J:1969|xxx:-||

\L|Abk:Martin|Sig:Ma|Au:Martin, P.|Tit:Temperley-Lieb algebras for
nonplanar statistical mechanics---the partition algebra
construction|Zs:J. Knot Theory Ramifications|Bd:3|S:51--82|J:1994|xxx:-||

\L|Abk:MS|Sig:MS|Au:Martin, P.; Saleur, H.|Tit:Algebras in
higher-dimensional statistical mechanics---the exceptional partition
(mean field) algebras|Zs:Lett. Math. Phys.|Bd:30|S:179--185|J:1994%
|xxx:hep-th/9302095||

\B|Abk:Saa|Sig:Sa|Au:Saavedra Rivano, N.|Tit:Cat\'egories
Tannakiennes|Reihe:Lecture Notes in Mathematics,
{\bf265}|Verlag:Springer-Verlag|Ort:Berlin-New York|J:1972|xxx:-||

\L|Abk:St1|Sig:St|Au:Stanley, R.|Tit:Modular elements of geometric
lattices|Zs:Algebra Universalis|Bd:1|S:214-217|J:1971/72|xxx:-||


\L|Abk:Wi|Sig:Wi|Au:Wilf, H.|Tit:Hadamard determinants, M\"obius
functions, and the chromatic number of a
graph|Zs:Bull. Amer. Math. Soc.|Bd:74|S:960--964|J:1968|xxx:-||

\endrefs

\bye